\documentclass[a4paper,12pt]{article}
 
	\usepackage[francais,english]{babel} 
	\usepackage[latin1]{inputenc} 
	\usepackage[T1]{fontenc} 
	\usepackage{epsfig}
	\usepackage{graphicx} 
	\usepackage{latexsym,array,multicol}
	\usepackage[T1]{fontenc}
	\usepackage[latin1]{inputenc}
	\usepackage{amsfonts}
	\usepackage{amsmath}
	\usepackage{amssymb,amsthm}
	\usepackage{textcomp}

\let\Bb\mathbb 
\let\g\mathfrak 
\let\Cl\mathcal 
\newcommand{\bi}{\begin{enumerate}}
\newcommand{\ei}{\end{enumerate}}
\newcommand{\vd}{\emptyset}
\newcommand{\bs}{\backslash}
\newcommand{\ra}{\rightarrow}
\newcommand{\fa}{\mapsto} 
\newcommand{\bd}{\partial}
\newcommand{\cx}{\times}
\newcommand{\grad}{{\rm grad}}
\newcommand{\gdiv}{{\rm div}}
\newcommand{\Vol}{{\rm Vol}}

\newcommand{\ch}{{\rm ch}}
\newcommand{\sh}{{\rm sh}}
\newcommand{\Th}{{\rm th}}
\newcommand{\inj}{{\rm inj}}
\newcommand{\beq}{\begin{equation}}
\newcommand{\eeq}{\end{equation}}
\newcommand{\pgh}{\paragraph{}} 
\newcommand{\Trace}{{\rm Trace}}
\newcommand{\inter}{\stackrel{\circ}} 

\newtheorem{theo}{Th\'eor\`eme}[section]
\newtheorem{lemm}[theo]{Lemme}
\newtheorem{prop}[theo]{Proposition}
\newtheorem{coro}[theo]{Corollaire}
\theoremstyle{definition}\newtheorem{defi}{D\'efinition}[section]
\theoremstyle{remark}\newtheorem{quest}{Question}
\theoremstyle{remark}\newtheorem{rema}{Remarque}[section]
\theoremstyle{remark}
\author{Samuel TAPIE}
\date{16/07/2008}

\title{Bas du spectre de surfaces hyperboliques de volume infini}

\begin{document}

\maketitle 

\selectlanguage{francais}

\begin{abstract}
Cet article présente des méthodes de contrôle du bas du spectre du Laplacien $\lambda_0$ sur des surfaces hyperboliques de volume infini. Nous commençons par donner une borne supérieure du $\lambda_0$ pour une surface géométriquement finie en fonction de la géométrie du coeur convexe. Nous nous intéressons ensuite à des surfaces de genre infini périodiques, construites en recollant des copies d'une surface géométriquement finie à bord selon le plan d'un graphe infini. Nous contrôlons le $\lambda_0$ de la surface infinie ainsi obtenue par des constantes issues des propriétés spectrales de la cellule élémentaire et des données combinatoires du graphe. Nous généralisons ensuite ces méthodes pour contrôler le $\lambda_0$ de deux autres types de surfaces de genre infini : celles qui admettent un découpage en morceaux bornés, et certains revêtements riemanniens.
\end{abstract}


\tableofcontents{}

\section*{Introduction}
Nous cherchons dans cette étude à relier certaines propriétés géométriques des surfaces de Riemann munies d'une métrique hyperbolique de volume total infini avec le bas du spectre du Laplacien associé à cette métrique. Nous emploierons indifféremment pour notre problème les termes \emph{surface de Riemann} et \emph{surface hyperbolique}, et cela sous-entendra toujours que nos variétés sont \emph{orientables}. Nous commençons par nous intéresser au cas  \emph{géométriquement fini}, c'est-à-dire lorsque le groupe fondamental de notre surface est de type fini. Elle se découpe alors canoniquement en une partie convexe de volume fini, son \emph{c\oe ur convexe}, à laquelle viennent se greffer des cylindres topologiques de volume infini que nous appellerons des \emph{vasques} (voir Section \ref{ssec:PreSurf}). Nous utilisons alors un résultat de \cite{Bu} pour contrôler la valeur du bas du spectre en fonction de l'aire (hyperbolique) du c\oe ur convexe et des longueurs des géodésiques qui le bordent, pour obtenir à la Section \ref{sec:GFin}

\begin{theo}\label{th:GFin}
Il existe une constante $R_2>0$ telle que si $M$ est une surface hyperbolique complète, non compacte et géométriquement finie, $C(M)$ son c\oe ur convexe et $\lambda_0(M)$ le bas de son spectre, alors
$$\lambda_0(M)\leq R_2 \frac{\ell (\bd C(M))}{\Vol (C(M))}.$$
\end{theo}

Nous utilisons pour cela une méthode identique à celle employée dans \cite{Ca} pour montrer un résultat analogue pour les variétés hyperboliques géométriquement finies de dimension 3. Nous montrons à la Section \ref{ssec:PincGeod} que lorsque l'on pince uniformément les géodésiques qui bordent le c\oe ur convexe, des méthodes développées par B. Colbois et Y. Colin de Verdière permettent d'obtenir un contrôle plus précis du bas du spectre en fonction des mêmes invariants. Ces deux méthodes nous permettent d'obtenir une condition suffisante pour l'existence d'une fonction propre associée au bas du spectre.

\pgh 

Nous nous intéressons ensuite à certaines surfaces hyperboliques dont le $\pi_1$ n'est pas de type fini, que nous appelons \emph{périodiques} au sens où elles reproduisent une infinité de fois la même \emph{cellule} selon le plan d'un \emph{graphe} à valence constante. Notons $\lambda_0(M)$ le bas du spectre d'une surface infinie périodique $M$, et $\lambda_0^N(C)$ le bas du spectre de la cellule élémentaire avec condition de Neumann aux bords. Toutes ces notions sont définies à la Section \ref{ssec:PreSpec}. Notre résultat principal est le suivant :

\begin{theo}\label{th:GInf1}
Soit $M$ une surface hyperbolique sans bord modelée sur un graphe $G$ à partir d'une cellule $C$ à trou spectral positif,
$$\lambda_0(M) \geq \lambda^N_0(C)$$
avec égalité si et seulement si $G$ est moyennable.
\end{theo}

On note $\g h(G)$ la constante de Cheeger du graphe, et $\mu_0(G)$ le bas du spectre de son Laplacien combinatoire (voir Section \ref{ssec:GInfMoy}). Notre démonstration donne en fait un résultat plus précis que celui-ci :

\begin{theo}\label{th:GInf2}
Soit $M$ une surface hyperbolique sans bord modelée sur un graphe $G$ à partir d'une cellule $C$ à trou spectral positif, on a 
$$\lambda_0^N(C)+A_1\mu_0(G)\leq \lambda_0(M) \leq \lambda_0^N(C)+A_2 \g h(G),$$
où  $A_1$ et $A_2>0$ sont des constantes qui ne dépendent que de la longueur des géodésiques de $\bd C$, du nombre de composantes de bord et de propriétés spectrales de la cellule.
\end{theo}

Nous montrons ensuite que la même méthode s'adapte et permet d'obtenir un résultat analogue pour des surfaces non périodiques qui admettent un \emph{découpage borné} (voir Section \ref{ssec:GInfBorn}) :

\begin{theo}\label{th:GInfBorn}
Soit $M$ une surface hyperbolique telle que 
$$M = \bigcup_{i\in \Bb Z} C_i,$$
où les $C_i$ sont des surfaces hyperboliques d'intérieurs disjoints à bords géodésiques, telles qu'il existe des constantes $k,K,\eta,v>0,$ avec
$$\forall i\in\Bb Z, k<\ell(\alpha)<K \mbox{ et } k<\Vol(C_i)<K,$$
où $\alpha$ parcours l'ensemble des composantes de $\bd C_i$, 
$$\forall i\in\Bb Z, \lambda_1^N(C_i)\geq\eta,$$
et le nombre de composantes de bord de $C_i$ est borné par $v$.

Il existe des constantes $A_1$ et $A_2$ ne dépendant que de $k,K,\eta$ et $v$ telles que
$$A_1\mu_0(G)\leq\lambda_0(M)\leq A_2\g h(G).$$
\end{theo}

Enfin, nous adaptons notre méthode pour l'appliquer au cas de certains revêtements riemanniens. Les travaux les plus connus reliant la moyennabilité d'un groupe de revêtement avec des résultats sur le spectre du Laplacien des variétés concernées sont certainement ceux de Robert Brooks  \cite{Br} et \cite{Br2}. Le résultat suivant est sans doute le plus abouti sur cette question, et le plus proche de nos considérations :

\begin{theo}[Brooks, 86]\label{th:Br}
Soit $p: M\ra M_1$ un revêtement riemannien galoisien, et $\Gamma = \pi_1(M_1)/\pi_1(M)$ son groupe de revêtement.
Si $M$ possède un domaine fondamental $F$ pour l'action de $\Gamma$ vérifiant la propriété \emph{(Br)}, alors $$\lambda_0(M) \geq \lambda_0(M_1)$$ avec égalité si et seulement $\Gamma$ est moyennable.
\end{theo}

Ce domaine $F$ joue ici le même rôle que notre cellule $C$ dans les énoncés précédent. La propriété \emph{(Br)} est une propriété technique sur les caractéristiques spectrales de $F$ difficile à contrôler, que nous explicitons à la Section \ref{ssec:GInfRev}. Les seuls exemples que donne Brooks de variétés vérifiant cette propriétés sont les variétés hyperboliques sans cusp à trou spectral positif. On retrouve donc des hypothèses très proches de celles de nos résultats. La Section \ref{ssec:GInfRev} détaille les recoupements et les différences de ces travaux avec les nôtres, qui aboutissent dans cette situation au théorème :

\begin{theo}\label{th:GInfRev}
Soit $M_1$ une surface hyperbolique à trou spectral positif, et $M\ra M_1$ un revêtement riemannien galoisien de groupe de revêtement $\Gamma$ de type fini. Supposons qu'il existe un domaine fondamental $F$ dans $M$ pour l'action de $\Gamma$ dont le bord est une union de géodésiques fermées.

Alors il existe des constantes $A_1$ et $A_2$ ne dépendant que de propriétés spectrales de $M_1$ et de la longueur des composantes de $\bd F$ telles que
$$\lambda_0(M_1)+A_1\mu_0(\Gamma)\leq \lambda_0(M)\leq \lambda_0(M_1)+A_2 \g h(\Gamma).$$
\end{theo}

Notre étude se conclut par la présentation de quelques questions que ce travail pose naturellement. 
\pgh
Nous présentons en appendice la démonstration du théorème clé qui permet de passer d'une surface construite par recollement de cellules à un revêtement riemannien. Il s'agit d'une caractérisation du bas du spectre du Laplacien sur une variété $M$ quelconque, avec condition de Neumann au bord, que nous notons $\lambda\leq\lambda_0^N(M)$ :

\begin{theo}[Sullivan, 87] \label{th:CarSul}
Pour tout réel $\lambda$, il existe une fonction $\phi$ $\Cl C^\infty$ \emph{$\lambda$-harmonique positive} sur $M$ avec condition de Neumann sur $\bd M$ si et seulement si $\lambda\leq\lambda_0^N(M)$.
\end{theo}

Ce résultat est dû à \cite{Sul} dans le cas sans bord, et nous reprenons avec plus de détails sa démonstration basée sur la théorie de la diffusion, en l'adaptant à la présence éventuelle d'un bord.
\pgh
Je remercie vivement Gilles et Gérard de m'avoir laissé vagabonder à la recherche de surfaces de genre infini sur lesquelles je pouvais dire quelque chose, ainsi que pour leur soutien au long de ce périple ; merci aussi à Constantin Vernicos et Françoise Dalbo pour leurs questions et remarques qui ont bien étoffé cette étude, et à Didier Piau pour son aide précieuse et patiente lors de la rédaction de l'appendice.

\section{Préliminaires}

Nous donnons ici quelques définitions et résultats élémentaires sur les surfaces hyperboliques et sur l'étude du spectre du Laplacien sur ces surfaces que nous emploierons par la suite. Bien qu'en principe suffisante pour comprendre les résultats de notre article, cette section ne prétend pas être une introduction complète à ces sujets, et il est vivement conseillé au lecteur intéressé de se reporter à la bibliographie citée.

\subsection{Surfaces hyperboliques}\label{ssec:PreSurf}

Une \emph{surface hyperbolique} $M$ est une variété de dimension 2 munie d'une métrique riemannienne à courbure constante égale à $-1$. Si $M$ est simplement connexe, alors $M$ est isométrique au disque unité muni de la métrique de Poincaré ; sinon $M$ est isométrique au quotient du disque par un sous groupe discret de $PSL_2(\Bb R)$. Toute surface hyperbolique compacte $M$ (éventuellement à bord géodésique) peut être réalisée comme réunion finie de \emph{pantalons hyperboliques} reliés deux à deux le long de bords géodésiques de longueurs identiques. La métrique hyperbolique de $M$ est complètement déterminée par les longueurs des bords de ces pantalons et un paramètre de raccordement pour chaque géodésique commune à deux pantalons, appelé \emph{angle de twist}. On pourra se reporter à \cite{BePe}, chapîtres B et D pour plus de détails. 

Une surface hyperbolique quelconque est la réunion de pantalons (compacts) dont les bords sont de longueurs positives, de pantalons non compacts dont certains bords sont de longueurs nulle et rejetés à l'infini (le voisinage d'un bord de longueur nulle, nécessairement non compact, est appelé un \emph{cusp}), et de composantes de volume infini que nous appellerons des \emph{vasques}, connues aussi sous le nom de \og \emph{funnel}\fg \ ou \og\emph{expanding ends}\fg, homéomorphes à des cylindres, dont un bord est une géodésique reliée à l'un des pantalons ou une autre vasque, tandis que l'autre bord est renvoyé à l'infini. 

\pgh

Nous aurons besoin par la suite d'écrire précisément la métrique hyperbolique sur un voisinage d'un cusp et au voisinage d'une géodésique fermée. Pour tout point $p$ d'une variété riemannienne $M$, on note $inj(p)$ le rayon d'injectivité de $M$ en $p$. On appelle partie $\epsilon$-\emph{mince} l'ensemble des points de $M$ où le rayon d'injectivité est inférieur à $\epsilon$, partie $\epsilon$-épaisse son complémentaire. Pour tout $\epsilon>0$, un cusp possède une partie $\epsilon$-mince ; voici une façon d'écrire sa métrique (voir \cite{BePe} p 151):

\begin{prop}
Soit $P_C$ un pantalon à cusp, $\epsilon>0$ inférieur à la longueur de la plus petite géodésique fermée de $P_C$, et $$V_\epsilon = \{p\in P_C : \inj(p)\leq \epsilon\}$$. Supposons $V_\epsilon$ connexe, c'est-à-dire que $P_C$ ne présente qu'un seul cusp. Alors $V_\epsilon$ est isométrique à $\Bb S^1\cx[0,\infty[$ muni de la métrique
$$ds^2 = e^{-2r}(\frac{\epsilon}{2\pi})^2d\theta + dr^2.$$
\end{prop}

\begin{coro}\label{coro:VolCusp}
Avec les notations précédentes, le volume de $V_\epsilon$ vaut $\epsilon$.
\end{coro}

Toute géodésique fermée dans une surface hyperbolique admet un voisinage homéomorphe à un cylindre que nous appellerons un collier, d'autant plus large que la longueur de la géodésique est petite. Sur ce voisinage, la métrique hyperbolique s'écrit sous la forme suivante (voir \cite{Cb}):

\begin{prop}[Lemme du Collier]\label{prop:LemCol}
Soit $\alpha$ une géodésique fermée de longueur $l$ contenue dans l'intérieur d'une surface hyperbolique $M$, alors $\alpha$ admet dans $M$ un voisinage isométrique à $\Bb S^1\cx]-m(l), m(l)[$, où $\alpha = \Bb S^1\cx\{0\}$, muni de la métrique
$$ds^2 = (\frac{l}{2\pi})^2\ch^2r d\theta + dr^2,$$
avec
$$m(l) = {\rm argsinh}(\frac{1}{\sh(l/2)}).$$
\end{prop}

Les coordonnées $(\theta, r)$ ainsi décrites s'appellent \emph{coordonnées de Fermi} sur le voisinage collier de $\alpha$. Si $\alpha$ est l'intersection entre un pantalon et une vasque, la métrique hyperbolique sur la vasque s'écrit également en coordonnées de Fermi 
$$ds^2 = (\frac{l}{2\pi})^2\ch^2r d\theta + dr^2,$$
pour $(\theta,r)\in\Bb S^1\cx [0,\infty[$. Si $\alpha$ est une composante de bord de $M$, $\alpha$ possède toujours un voisinage (relatif) collier dans $M$ qui s'écrit $(\theta,r)\in\Bb S^1\cx[0,m(l)[$ muni de la métrique de Fermi.

\pgh

\begin{defi}
Soit $M$ une variété hyperbolique (de dimension quelconque), son \emph{c\oe ur convexe} $C(M)$ est le plus petit convexe de $M$ tel que $C(M)$ soit homéomorphe à $M$ (voir Fig.1). On dit que $M$ est \emph{géométriquement finie} si et seulement si le volume de son coeur convexe est fini.
\end{defi}

\begin{figure}
\begin{center}
\includegraphics[width = 0.6\textwidth]{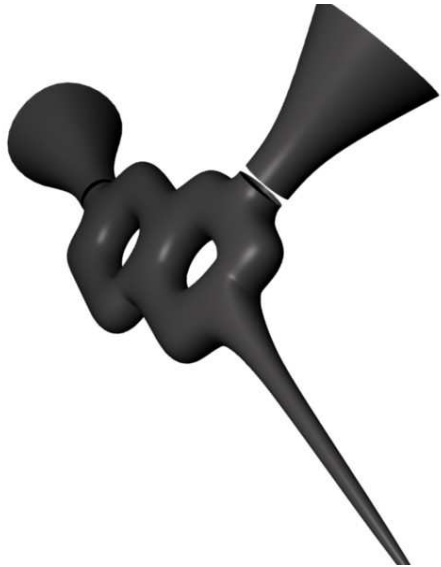}
\caption{La composante connexe à l'avant plan est le c\oe ur convexe de la surface, avec un cusp}
\end{center}
\end{figure}

Dans le cas des surfaces hyperboliques qui nous intéresse, le c\oe ur convexe n'est rien d'autre que $M$ privée de toutes ses vasques. Comme le volume d'un pantalon, éventuellement à cusp, est fini uniformément minoré (c'est un invariant topologique par la formule de Gauss-Bonnet), nous voyons que le volume de $C(M)$ est fini si et seulement si $M$ est la réunion d'un nombre fini de pantalons, pantalons à cusps et vasques, c'est à dire lorsque \emph{son groupe fondamental est de type fini}. Ces deux notions ne sont plus équivalentes en dimension plus grande, voir par exemple \cite{Ca}.

Nous travaillerons également sur des surfaces qui ne sont pas géométriquement finies. Nous supposerons toujours qu'elles sont réunion dénombrable de compacts : elles sont alors nécessairement réunion dénombrable de pantalons compacts, de pantalons à cusps et de vasques. 

Dans le cas géométriquement fini, nous nous efforcerons à la Section \ref{sec:GFin} de relier le volume et la longueur du bord du c\oe ur convexe aux éléments spectraux que nous présentons maintenant.

\subsection{Eléments de théorie spectrale}\label{ssec:PreSpec}

Le lecteur voulant entrer dans les détails de ce que nous présentons maintenant sans démonstration pourra se référer à \cite{Cha}.

Soit $M$ une variété riemannienne. Nous appelons \emph{Laplacien} l'opérateur $\Delta$ défini sur toute fonction $\Cl C^2$ sur $M$ (à valeur réelle) par 
$$\Delta f = \gdiv(\grad f) = -\Trace(\nabla_. (\nabla f)),$$
où $\nabla$ est la connexion de Levi-Civita (et donc le gradient usuel lorsqu'elle s'applique à une fonction différentiable sur $M$). Notons que nous choisissons pour le Laplacien une convention de signe opposée à celle utilisée dans \cite{Cha} et dans la plupart des publication américaines. Notre convention, utilisée par nombre de géomètres français, a l'avantage de donner un opérateur défini positif, comme le montre la Formule de Green ci-dessous. Le lecteur est de toutes façons invité à se méfier fortement des signes toutes les fois qu'il se réfèrera à un article traitant du Laplacien, quelle que soit la nationalité de son auteur !!

A l'aide de l'expression des opérateurs gradients et divergence, on peut exprimer localement le Laplacien comme opérateur différentiel dans un système de coordonnées à partir de l'expression de la métrique. La seule expression explicite qui nous intéssera ici est celle du Laplacien sur le voisinage collier d'une géodésique fermée de longueur $l$ en coordonnées de Fermi : la formule de \cite{Cha}, p 5 devient
\beq \label{eq:LapFermi} \Delta f = -(\frac{\bd^2f}{\bd r^2}+\Th r \frac{\bd f}{\bd r} + \frac{(2\pi)^2}{l^2\ch^2r}\frac{\bd^2f}{\bd\theta^2}).\eeq

Si $M$ est une variété à bord, nous dirons qu'une fonction $f$ de classe $\Cl C^1$ sur $M$ vérifie les \emph{conditions de Dirichlet} si elle est nulle sur $\bd M$, et les \emph{conditions de Neumann} si son gradient vérifie en tout point $x\in\bd M$ :
$$g_x(\nabla f(x),\nu_x) = 0$$
où $\nu_x$ est la normale au bord au point $x$ considéré. Sauf dans notre appendice, nous nous placerons toujours dans l'un de ces deux cas, ce qui nous permet d'écrire le résultat suivant sans terme de bord.

\begin{prop}[Formule de Green]\label{prop:Green}
Soient $f,g$ deux fonctions de classe $\Cl C^2$ sur $M$ vérifiant les conditions de Neumann ou de Dirichlet sur $\bd M$ si $\bd M\neq \vd$, alors
$$\int_M\Delta f g = \int_M \nabla f.\nabla g.$$
\end{prop}

On définit $\Cl H^1(M)$ comme l'ensemble des fonctions $f\in L^2(M)$ telles que le gradient de $f$ au sens des distributions est un champ de vecteur de $L^2(M)$ que nous noterons encore $\nabla f$.

$\Cl H^1(M)$ muni de la norme
$$||f||^2_{\Cl H^1(M)} = ||f||^2_{L^2(M)} + ||\nabla f||^2_{L^2(M)}$$
est alors un espace de Hilbert, et l'ensemble des fonctions $\Cl C^\infty$ à support compact dans $M$ est dense dans $\Cl H^1(M)$. On définit l'espace $\Cl H_0^1(M)$ comme le complété dans $\Cl H^1(M)$ de l'ensemble des fonctions $\Cl C^\infty$ à support compact dans $\stackrel{\circ}{M}$.  On peut alors définir sur $\Cl H^1(M)$ ou sur $\Cl H^1_0(M)$ la forme quadratique \emph{énergie de Dirichlet} (ou simplement \emph{énergie}) :
$${\rm Dir}(f) = ||\nabla f||^2_{L^2(M)}.$$
L'opérateur non compact qui est associé à cette forme quadratique est l'extension de Friedrich du Laplacien (défini précédemment pour des fonctions $\Cl C^2$) à $\Cl H^1(M)$. Lorsque $M$ a un bord, si nous effectuons cette extension  sur  $\Cl H^1(M)$, nous obtenons le Laplacien avec conditions de Neumann, tandis que si nous la limitons à  $\Cl H^1_0(M)$, nous obtenons le Laplacien avec conditions de Dirichlet. Nous noterons dans tous les cas $\Delta$ ces extensions et le contexte précisera si nous travaillons avec des fonctions $\Cl C^2$ ou dans $\Cl H^1$, ainsi que les conditions au bord considérées. 

\pgh

Le \emph{spectre} du Laplacien est l'ensemble des $\lambda\in\Bb R$ tels que $\Delta-\lambda$ vu comme opérateur sur un espace hilbertien ($\Cl H^1(M)$ ou $\Cl H^1_0(M)$ suivant les situations) n'est pas inversible. D'après la formule de Green, c'est un sous-ensemble de $\Bb R^+$. On appelle \emph{bas du spectre} sa borne inférieure, notée $\lambda_0(M)$ pour $M$ sans bord, $\lambda^N_0(M)$ pour le bas du spectre avec condition de Neumann, et $\lambda^D_0(M)$ pour le bas du spectre avec condition de Dirichlet. On a toujours :
si $\bd M = \vd$,
$$\lambda_0(M) = \inf_f\{\frac{||\nabla f||^2}{||f||^2}\}$$
où $f$ parcourt l'ensemble des fonctions de $\Cl H^1(M)$ ; si $\bd M\neq \vd$, 
$$\lambda_0^N(M) = \inf_f\{\frac{||\nabla f||^2}{||f||^2}\}$$
où $f$ parcourt l'ensemble des fonctions de $\Cl H^1(M)$, et
$$\lambda_0^D(M) = \inf_f\{\frac{||\nabla f||^2}{||f||^2}\}$$
où $f$ parcourt l'ensemble des fonctions de $\Cl H^1_0(M)$. Pour toute fonction $f$ de $\Cl H^1$ (par exemple $\Cl C^1$ par morceaux), on appelle $$\frac{||\nabla f||^2}{||f||^2}$$ son \emph{quotient de Rayleigh}.

On dit qu'une fonction $f$ est $\lambda$\emph{-harmonique} si $\Delta f = \lambda f$, et \emph{fonction propre} du Laplacien associée à la valeur propre $\lambda$ si elle est dans $\Cl H^1$ et $\lambda$-harmonique ; $\lambda$ est alors appelée \emph{valeur propre} du Laplacien. Une valeur propre est nécessairement un point du spectre, donc positive. Si $M$ est \emph{compacte}, on montre que le spectre est l'ensemble (discret) de la suite de ses valeurs propres, qui sont alors de multiplicité finie. Pour $M$ de \emph{volume fini}, $\lambda_0 = 0$ est valeur propre associée aux fonctions constantes. Lorsque $M$ n'est pas de volume fini, l'existence de fonctions propres (et donc de valeurs propres) n'est pas assurée. Le résultat suivant, que nous utiliserons souvent par la suite, regroupe plusieurs théorèmes classiques :

\begin{theo}\label{th:EigFunc}
S'il existe une fonction $\psi_0\in \Cl H^1(M)$ (resp. dans $\Cl H_0^1(M)$) telle que son quotient de Rayleigh soit égal à $\lambda_0^N(M)$ (resp. $\lambda_0^D(M)$).

Alors $\psi_0$ est strictement positive sur $\inter{M}$, de classe $\Cl C^\infty$ sur $M$, et est fonction propre du Laplacien avec condition de Neumann (resp. de Dirichlet) aux bords. Toute fonction propre du Laplacien associée à la valeur $\lambda_0$ est donc proportionnelle à $\psi_0$.
\end{theo}

\begin{defi}\label{def:Trou}
Soit $M$ une variété vérifiant les hypothèses du théorème 1.6, on note
$$\lambda_1 = \inf\{\frac{||\nabla g||^2}{||f||^2}\int_M\psi_0 g = 0\}$$
où les fonctions $f$ sont prises dans $\Cl H^1(M)$ ou $\Cl H^1_0(M)$ suivant le cas considéré, et on appelle $\eta = \lambda_1-\lambda_0\geq 0$ le \emph{trou spectral} de $M$.
\end{defi}

Lorsque $M$ est compacte (avec ou sans bord) le trou spectral de $M$ est donc strictement positif. Dans la suite de notre étude, nous découperons régulièrement les surfaces étudiées en morceaux disjoints. Le lemme élémentaire suivant nous sera alors utile :

\begin{lemm}\label{lemm:MonoNeum}
Si $M'\subset M$ sont deux variétés complètes à bords compacts, alors $$\lambda_0^N(M)\geq\min(\lambda_0^N(M'),\lambda_0^N(M\bs M')).$$
\end{lemm}

\begin{proof}
Soit $f\in \Cl H^1(M)$, on a  $\hat f = f_{|M'}\in \Cl H^1(M')$, et son gradient (au sens des distributions) sur $M'$ est un champ de vecteurs qui vérifie $$\nabla \hat f(x) = \nabla f(x).$$ De même, $f_{|M\bs M'}\in \Cl H(M\bs M')$.
De plus, pour tous nombres positifs $a,b,c,d$,
$$\frac{a+b}{c+d}\geq\min\{\frac{a}{c},\frac{b}{d}\}.$$
On a donc
$$\frac{||\nabla f||^2}{||f||^2}\geq \min\{\frac{||\nabla \hat f||^2}{||\hat f||^2}, \frac{||\nabla (f-\hat f)||^2}{||f-\hat f||^2}\}\geq\min\{\lambda_0^N(M'),\lambda_0^N(M\bs M')\}.$$
Ce résultat étant valable pour toute fonction $f\in \Cl H^1(M)$, on a bien 
$$\lambda_0^N(M)\geq\min\{\lambda_0^N(M'\},\lambda_0^N(M\bs M')).$$
\end{proof}

\pgh
Certaines de nos surfaces auront des propriétés de symétrie, que nous utiliserons pour nos problèmes spectraux à l'aide du résultat suivant :

\begin{prop}\label{prop:InvEigFunc}
Soit $M$ une variété riemannienne (éventuellement à bord) munie d'une isométrie $I$ d'ordre fini $v$. Si le bas du spectre de $M$ est atteint par une fonction $\psi_0$, alors $\psi_0$ est invariante par $I$.
\end{prop}

\begin{proof}
D'après le Théorème \ref{th:EigFunc}, $\psi_0$ est l'unique fonction propre associée à la valeur propre $\lambda_0$. Or, posons
$$\Psi_0 = \psi_0 + \psi_0\circ I+...+\psi_0\circ I^{v-1},$$ cette fonction est invariante par $I$ et on a encore
$$\Delta \Psi_0 = \lambda_0 \Psi_0.$$
Il existe donc une constante $k$ telle que $\Psi_0 = k\psi_0$, qui ne peut être nulle puisque $\psi_0$ est strictement positive sur $\inter{M}$ d'après le Théorème \ref{th:EigFunc}. On a donc $\psi_0\circ I = \psi_0$.
\end{proof}

\begin{coro}\label{coro:InvEigFunc}
Soit $M$ une variété riemannienne (éventuellement à bord avec condition de Neumann) munie d'une isométrie $I$ d'ordre fini $v$,
alors 
$$\lambda_0(M) = \inf_f\{\frac{||\nabla f||^2}{||f||^2}\},$$
où $f$ parcourt l'ensemble des fonctions $\Cl C^\infty$ à support compact dans $M$ invariantes par $I$.
\end{coro}

\begin{proof}
Soit $(U_i)_i$ une suite croissante d'ouverts à fermetures compactes de $M$ tels que $$M = \bigcup_i U_i,$$
il existe une suite d'ouverts à fermetures compactes $(V_i)_i$ invariants par $I$ tels que pour tout $i$,
$$U_i\subset V_i.$$
En particulier, $M$ est l'union des $\{V_i\}$. Comme l'ensemble des fonctions à support compact dans $M$ est dense dans $\Cl H^1(M)$, on a 
$$\lambda_0(M) = \inf_i\inf_{f\in\Cl C^\infty V_i}\frac{||\nabla f||^2}{||f||^2}.$$
D'après la proposition précédente, pour tout $i$,
$$\inf_{f\in\Cl C^\infty V_i}\frac{||\nabla f||^2}{||f||^2}$$
est atteint par une fonction invariante par $I$, ce qui conclut notre démonstration.
\end{proof}

\subsection{Géométrie spectrale des surfaces hyperboliques}

Par la suite, nous ne travaillerons que sur des surfaces hyperboliques. Nous présentons ici quelques uns des résultats connus sur le spectre du Laplacien sur ces surfaces, qui sont à l'origine des motivations de cet article, et permettent en particulier de comprendre les hypothèses de nos théorèmes dans un cadre plus général de géométrie hyperbolique.

\begin{rema}
La plupart des résultats présentés dans cette section sont valables en dimension plus grande en adaptant simplement certaines constantes. Pour plus de détails à ce sujet, le lecteur est invité à consulter les références citées.
\end{rema}

Le résultat suivant démontré dans \cite{Do} assure, suivant la valeur de $\lambda_0$, l'existence d'une fonction propre associée sur une surface de Riemann :

\begin{theo}\label{th:Don}
Soit $M$ une surface hyperbolique géométriquement finie non compacte, alors la demi-droite $[1/4,\infty[$ est dans le spectre du Laplacien, et tout point du spectre dans $[0,1/4[$ (s'il y en a) est associé à une valeur propre de multiplicité finie.
\end{theo}

\begin{coro}\label{coro:Don}
Sous les hypothèses du théorème précédent, on a :
$$\lambda_0(M)\leq 1/4 ;$$
Si $\lambda_0(M)<1/4$, alors $\lambda_0$ est une valeur propre simple associé à une fonction propre positive $\psi_0$ et $$\lambda_0<\lambda_1\leq\frac{1}{4} :$$
le trou spectral de $M$ est donc strictement positif et inférieur à $1/4$. 
\end{coro}

Dans cet énoncé, le $\lambda_1$ et le trou spectral correspondent à la Définition \ref{def:Trou}. Dans le cas géométriquement fini, on a donc une distinction importante entre la situation $\lambda_0<1/4$ et $\lambda_0 = 1/4$. Cette distinction peut aussi se comprendre à l'aide du résultat suivant, valable pour $M$ géométriquement finie ou infinie, démontré par exemple dans \cite{Sul} :

\begin{theo}[Sullivan]
Soit $M = \Bb H^2/\Gamma$ et $\delta$ l'exposant critique de Poincaré de $\Gamma$, on a
$$\lambda_0(M) = \left\{\begin{array}{ccc}
\delta(1-\delta) & \mbox{si} & \delta>1/2\\
1/4 & \mbox{si} & \delta \leq 1/2.\end{array}\right.$$
\end{theo}

\begin{coro}
On a toujours $\lambda_0(M) \leq 1/4$.
\end{coro}

Si $M$ est géométriquement finie, $\delta$ est aussi la \emph{dimension de Hausdorff de son ensemble limite}. Ces notions d'exposant critique de Poincaré et d'ensemble limite d'une surface hyperbolique sont bien connues des personnes étudiant la géométrie hyperbolique et la \emph{théorie ergodique des groupes discrets}. Le théorème de Sullivan permet d'exprimer la plupart de nos résultats sur le bas du spectre comme des résultats sur l'exposant critique, mais nous n'utiliserons pas ce point de vue dans cette étude. Le lecteur intéressé par ces question est donc invité à se référer, par exemple, à \cite{Sul} ou \cite{Ca}.

A la Section \ref{sec:GInf}, nous nous intéresserons à des surfaces à bord géodésique dont le trou spectral est positif. Les résultats ci-dessus s'appliquent dans le cas à bord : il suffit de faire le double de la surface le long de son bord (voir Corollaire \ref{coro:L0Pos}). Nos exemples principaux seront donc des surfaces compactes ou des surfaces géométriquement finies non compactes dont le bas du spectre est plus petit que $1/4$. Un calcul explicite de l'exposant critique montre que pour toute surface hyperbolique $M$, dès que $M$ présente un cusp $\delta>1/2$ et donc $\lambda_0(M)<1/4$ lorsque $M$ est géométriquement finie (voir par exemple \cite{McM} p 7). En l'absence de cusp, le Théorème \ref{th:GFin} que nous démontrons maintenant donne une condition suffisante pour avoir $\lambda_0<1/4$.

\section{Cas géométriquement fini}\label{sec:GFin}
\subsection{Contrôle du bas du spectre par la géométrie du c\oe ur convexe}

L'objectif de ce paragraphe est de démontrer le Théorème \ref{th:GFin} qui permet, dans le cas général d'une surface hyperbolique non compacte géométriquement finie, de contrôler le bas du spectre à partir de la géométrie de son c\oe ur convexe.  Ce paragraphe est une adaptation aux surfaces hyperboliques du résultat démontré par R.D.Canary dans \cite{Ca} pour les variétés hyperboliques de dimension 3 : bien qu'un analogue du théorème \ref{th:LemBus} soit déjà contenu dans [Buser], la démonstration que nous présentons ici est due à Canary. Rappelons que nous utilisons une convention de signe pour le Laplacien différente de celle de Canary.

\begin{lemm}\label{lemm:PreBus}
Pour tout $n\in \Bb N$, $n\geq 2$, il existe une constante $R_n$ t.q. pour toute variété complète non compacte $M$ dont la courbure de Ricci est minorée par $-(n-1)$, si $f$ est une fonction $\lambda$-harmonique sur $M$, $\lambda>0$, alors 
$$\frac{|\nabla f|^2}{f^2}\leq R_n^2.$$
\end{lemm}

\begin{proof}
D'après le théorème 1.2 de \cite{LY} rappelé dans \cite{Ca}, si $u(x,t)$ est une solution positive de l'équation de la chaleur $(\Delta + \frac{\bd}{\bd t})u(x,t) = 0$ sur $M\cx (0,\infty)$ où $M$ est une variété de dimension $n$ sans bord dont la courbure de Ricci est bornée inférieurement par $-(n-1)$, alors pour tout $\alpha>1$, 
$$\frac{|\nabla u|^2}{u^2} - \frac {\alpha u_t}{u} \leq \frac{n\alpha^2(n-1)}{\sqrt{2}(\alpha-1)} + \frac{n\alpha^2}{2t}.$$

Si f est une fonction $\lambda$-harmonique positive sur $M$, $u(x,t) = e^{-\lambda t}f(x)$ est une solution positive de l'équation de la chaleur. En posant $\alpha = 2$ et en faisant tendre $t$ vers l'infini, la majoration précédente donne alors : 
$$\frac{|\nabla f|^2}{f^2}\leq 2\sqrt{2}n(n-1) - 2\lambda.$$
Il suffit donc de poser 
$$R_n^2 = 2\sqrt{2}n(n-1).$$
\end{proof}

\begin{theo}[Lemme de Buser]\label{th:LemBus}
Si $M$ est une variété de dimension $n$ complète non compacte sans bord dont la courbure de Ricci est bornée inférieurement par $-(n-1)$, alors 
$$\lambda_0(M)\leq R_n \g h(M)$$
où $\g h(M)$ est la constante de Cheeger de $M$ et $R_n$ ne dépend que de $n$.
\end{theo}

\begin{proof}
On note $V$ la mesure riemannienne sur $M$ et $A$ la mesure induite sur les sous-variétés de $M$ de co-dimension $1$ et on rappelle qu'alors, la \emph{Constante de Cheeger de} $M$ est définie comme
$$\g h(M) = \inf_{M'}\frac{A(\bd M')}{V(M')},$$
où $M'$ parcourt l'ensemble des domaines compacts de $M$. On ne perd pas en généralité à supposer que $\lambda_0>0$. On sait d'après le Théorème \ref{th:CarSul} qu'il existe alors une  fonction $\lambda_0$-harmonique positive $f$ sur $M$. D'après le Lemme \ref{lemm:PreBus}, il existe une constante $R_n$ ne dépendant que de $n$ tq 
$$|\frac{\nabla f}{f}(x)|\leq R_n.$$
De plus, comme $\nabla(\log f) = \nabla f/f$, on a
$$\Delta \log f = \lambda_0 + \frac{|\nabla f|^2}{f^2}\geq \lambda_0.$$

Soit $M'$ un domaine de $M$ relativement compact à bord $\Cl C^1$ par morceaux, d'après la formule de Stokes,
$$\int_{M'}\Delta(\log f)dV  = \int_{\bd M'} \frac{1}{f}\frac{\bd f}{\bd \nu_w}dA(w),$$
avec par définition
$$\frac{1}{f}\frac{\bd f}{\bd \nu_w} = g_w(\frac{\nabla f}{f},\nu_w) :$$
on a donc 
$$\int_{\bd M'} \frac{1}{f}\frac{\bd f}{\bd \nu_w}dA(w) \leq R_n A(\bd M').$$
Comme de plus, d'après ci-dessus,
$$\int_{M'}\Delta(\log f)dV\geq \lambda_0 V(M'),$$
on obtient que pour tout domaine $M'$ de $M$ à bord $\Cl C^1$ par morceaux
$$\lambda_0\leq R_n\frac{A(\bd M')}{V(M')},$$
d'où par définition de $\g h(M)$ : 
$$\lambda_0\leq R_n\g h(M).$$

\end{proof}

Le théorème \ref{th:GFin} s'énonce alors précisément :

\begin{theo}
Si $M$ est une surface hyperbolique complète, non compacte et géométriquement finie, on a
$$\lambda_0(M)\leq R_2 \frac{\ell (\bd C(M))}{\Vol (C(M))}.$$
\end{theo}
\begin{proof}
Pour $\epsilon>0$, posons $$C_\epsilon(M) = C(M)\cap M_e(\epsilon)$$ où $M_e(\epsilon)$ désigne la partie $\epsilon$-épaisse de $M$ (voir Section \ref{ssec:PreSurf}). Pour $\epsilon$ suffisamment petit (en particulier plus petit que la longueur de la plus petite géodésique fermée), $C_\epsilon(M)$ est le c\oe ur convexe de $M$ dont on a retiré la partie $\epsilon$-mince de ses cusps , qui sont en nombre fini $N_0$. On a donc 
$$\ell(\bd C_\epsilon(M)) = \ell(\bd C(M))+N_0\epsilon$$
et d'après le Corollaire \ref{coro:VolCusp},
$$\Vol(C_\epsilon(M)) = \Vol(C(M)) - N_0\epsilon.$$
Comme $C_\epsilon(M)$ est compact, on a pour tout $\epsilon>0$, 
$$\g h(M)\leq \frac{\ell(\bd C_\epsilon(M))}{\Vol(C_\epsilon(M))} :$$
le résultat désiré suit en faisant tendre $\epsilon$ vers $0$.
\end{proof}

\subsection{Pincements de géodésiques et spectre de graphes}\label{ssec:PincGeod}
On présente ici très succinctement les résultat de la méthode développée par Bruno Colbois et Yves Colin de Verdière dans \cite{Cb}, puis \cite{CbCV}. 

Soit $M$ une surface hyperbolique complète sans bord, $\gamma_1,...,\gamma_k$ des géodésiques fermées de $M$ découpant $M$ en N+n morceaux $M_1,...,M_N, M_{N+1},...,M_{N+n}$, où on suppose que pour $1\leq i\leq N$, les $M_i$,  sont de volume fini et pour $i>N$, les $M_i$ sont de volume infini. Soit $G$ le graphe à $N+n$ sommets, où deux sommets $i$ et $j$, $1\leq i<j\leq N+n$ sont reliés si et seulement si $M_i\cap M_j\neq\vd$. On a alors $M_i\cap M_j = \gamma_\alpha$, avec $1\leq\alpha\leq k$, et on assigne à l'arête $\{i,j\}$ la longueur $l_\alpha = \ell(\gamma_\alpha)$. On a donc un graphe à $N+n$ sommets et $k$ arêtes.

On note l'ensemble des sommets $S = S_1\cup S_2$, où $S_1 = \{1,...,N\}$ correspond à ceux associés aux parties de volume fini, et $S_2$ aux autres. On considère l'ensemble des fonctions de $S$ dans $\Bb R$, nulles sur $S_2$, de carré sommable pour la mesure $$\mu = \sum_{i\in{S}} V_i \delta(i)$$ où $V_i = \Vol(M_i)$ si $i\in S_1$, $V_i = 0$ sinon ; on note cet ensemble $L^2(S_1,\mu)$. 

Nous appellons dans cette partie \emph{Laplacien combinatoire} sur $G$ l'opérateur linéaire associé à la forme quadratique $q$ sur $L^2(S_1,\mu)$ définie par
$$q(f) = \frac{1}{\pi}\sum_{\alpha = 1}^k l_\alpha (f(i_\alpha)-f(j_\alpha))^2,$$
où $i_\alpha$ et $j_\alpha$ sont les extrémités de l'arête correspondant à la géodésique de découpage $\alpha$. En supposant $S_2$ non vide, cet opérateur a $N$ valeurs propres
$$0\leq\mu_0\leq \mu_1...\leq\mu_{N-1}.$$

Soit $\epsilon>0$ et $M^\epsilon$ l'unique surface hyperbolique obtenue en imposant à chaque géodésique $\gamma_\alpha$ la longueur $\epsilon l_\alpha$ (voir Section \ref{ssec:PreSurf}). 

\begin{theo}[Colbois, Théorème 1]\label{th:EqCol}
En gardant les notations précédentes, pour $\epsilon$ suffisamment petit, les $N$ premières valeurs propres simples du Laplacien sur $M^\epsilon$ (éventuellement égales, associées à des vecteurs propres indépendants) $\lambda^\epsilon_0\leq...\leq \lambda^\epsilon_{N-1}$ existent et vérifient lorsque $\epsilon\ra 0$
$$\lambda^\epsilon_i\sim \epsilon \mu_i.$$
\end{theo}

Si on applique ce résultat en pinçant uniformément toutes les géodésique qui forment le bord du c\oe ur convexe, on obtient à la limite un résultat plus précis que le Théorème \ref{th:GFin} :

\begin{prop}\label{prop:GFinPinc}
$$\lambda_0^\epsilon \sim \epsilon \frac{\ell(\bd C(M))}{\pi \Vol(C(M))}.$$
\end{prop}

\begin{proof}
Ici, $\cup_{\alpha} \gamma_\alpha = \bd C(M)$ ; le c\oe ur convexe est le seul morceau de volume fini, les autres (des vasques hyperboliques) sont de volume infini. On a donc $S_1 = \{s\}$, et pour toute fonction $x : S_1\ra \Bb R$, on a
$$q(x) = \frac{1}{\pi}\sum_\alpha l_\alpha (x(s))^2 = \frac{1}{\pi}\ell(\bd C(M))x^2$$
La première valeur propre de $q$ pour le produit scalaire 
$$\mu(x,y) = \sum_{i\in S_1} V_i x(i)y(i) = \Vol(C(M))xy$$
est donc $$\mu_0 = \frac{\ell(\bd C(M))}{\pi \Vol(C(M))},$$ le Théorème \ref{th:EqCol} nous donne alors 
$$\lambda_0^\epsilon \sim \epsilon \frac{\ell(\bd C(M))}{\pi \Vol(C(M))}.$$
\end{proof}

\begin{rema}
Les résultats que nous établirons dans la section suivante nécessiterons des surfaces (à bords) dont le trou spectral (avec condition de Neumann) est positif. Le Théorème \ref{th:GFin} et la Proposition \ref{prop:GFinPinc} assurent que $\lambda_0<1/4$ dès que la longueur du bord du c\oe ur convexe devient petite devant son volume. D'après le Théorème \ref{th:Don}, c'est une condition suffisante pour qu'il existe une fonction propre associée au bas du spectre et que le trou spectral soit positif.
\end{rema}

\section{Surfaces géométriquement infinies périodiques}\label{sec:GInf}

Nous présentons dans cette section des résultats de contrôle du spectre de certaines surfaces géométriquement infinie, c'est à dire de genre infini (voir Section \ref{ssec:PreSurf}). A partir d'une surface hyperbolique géométriquement finie (la \emph{cellule}) et d'un graphe, nous allons contruire toute une famille de surfaces hyperboliques géométriquement infinies, dont certaines auront un bas du spectre strictement positif : nous le caractériserons alors en fonction du bas du spectre de la cellule et du graphe sous-jacent.

\subsection{Construction d'une surface modelée sur un graphe à valence constante}

Soit $C$ une surface hyperbolique, pas forcément compacte ni géométriquement finie,  dont le bord est constitué de $v$ géodésiques fermées distinctes et de même longueur : $$\bd C = \alpha^1\cup \cdots \cup \alpha^v,$$ où les $\alpha^p$ sont disjointes. On suppose que $C$ est muni d'une isométrie $J$ d'ordre $v$ telle que pour tout $1\leq p\leq v-1$, 
$$J(\alpha^p) = \alpha^{p+1},$$ et que le trou spectral
$$\eta = \lambda_1^N(C)-\lambda_0^N(C)$$ est positif (voir Section \ref{ssec:PreSpec}). Cela signifie en particulier que $C$ admet une fonction propre $\psi_0$ associée à la valeur $\lambda_0$. C'est toujours le cas lorsque $M$ est de volume fini (compacte ou avec cusp d'après la section \ref{ssec:PreSurf}) ; dans le cas où $C$ est géométriquement finie de volume infini, si la longueur du bord du c\oe ur convexe est suffisamment petite devant son volume, la remarque qui conclut la section précédente donne des conditions pour que ce soit le cas.

On considère un graphe $G = (V,E)$ de valence constante $v$, dont les sommets sont indexés par $\Bb N$ : $V = \{x_i\}_{i\in\Bb N}$, et on note $i\sim j$ si et seulement si $(i,j)\in E$. On se donne une famille de copies de $C$ notées $\{C_i\}_{i\in\Bb N}$, dont on note les composantes de bord $(\alpha^k_i),_{1\leq k\leq v}$, et une famille d'isométries $\phi_i : C_i\ra C$ telles que pour tous $1\leq p\leq v$, 
$$\phi_i(\alpha_i^p) = \alpha^p$$.

\begin{defi}\label{def:SInfPer}
On dira qu'une surface hyperbolique $M$ est \emph{modelée sur le graphe $G$ à partir de la cellule $C$} si et seulement si on a 
$$M = \coprod_{i\in\Bb N} C_i/\sim,$$
où $\sim$ est une relation d'équivalence qui identifie $\alpha_i^p$ à $\alpha_{j}^q$ via $(\phi_j)^{-1}\circ J^{q-p}\circ \phi_i$ si et seulement si $i\sim j$ et $\alpha_i^p$ et $\alpha_{j}^p$ ne sont identifiés à aucun autre $\alpha_k^r$ (voir Fig. 2 $\&$ 3).
\end{defi}

\begin{rema}
Nous n'affirmons pas que, si la cellule et le graphe sont fixés, notre construction définit une surface de genre infini de façon unique (ce qui nécessiterait de définir ce que l'on entend par \emph{unique}). Pour la suite de notre étude, il nous suffit de savoir que notre surface de genre infini peut être obtenue de cette manière.
\end{rema}

\begin{rema}
On pourrait supposer que chacune des $v$ composantes $\alpha_i$ de la cellule $C$ est composée de plusieurs composantes connexes, chacune géodésique, et que l'isométrie $J$ échange les $\alpha_i$ globalement. Tous nos résultats restent valables dans ce cas : les démonstrations s'adaptent aisément, mais elles prennent une lourdeur inutile. Pour plus de clarté, nous supposerons toujours par la suite que chaque composante de bord $\alpha_i$ est composée d'une unique géodésique.
\end{rema}

\begin{figure}
\begin{center}
\includegraphics[width = 0.6\textwidth]{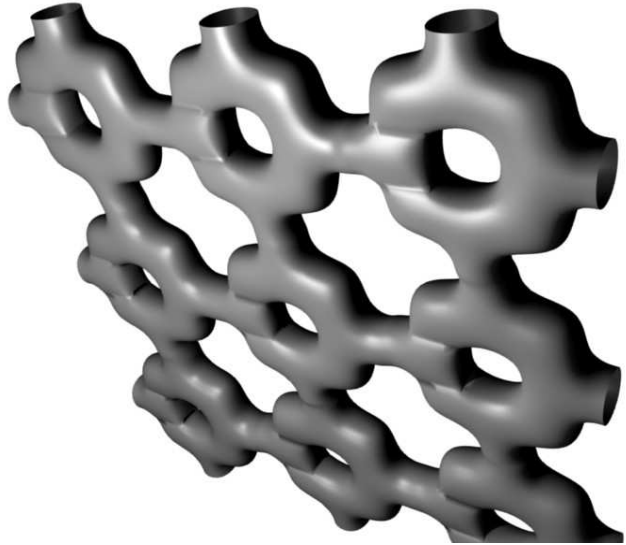}
\caption{Surface modelée sur $\Bb Z^2$ à partir d'un tore à $4$ composantes de bord}
\end{center}
\end{figure}

\begin{figure}
\begin{center}
\includegraphics[width = 0.6\textwidth]{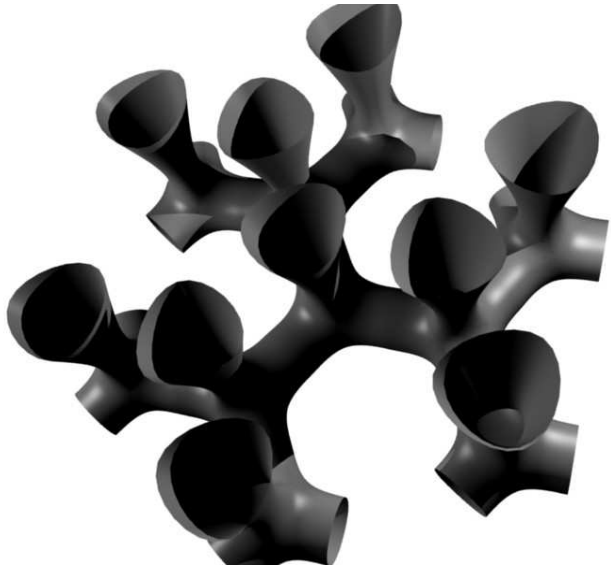}
\caption{Surface modelée sur l'arbre de degré 3 à partir d'une cellule de volume infini}
\end{center}
\end{figure}

On notera toujours $C_i = C_i/\sim$ chacune des parties de $M$, toutes isométriques à $C$. Nous noterons désormais le bas du spectre de la cellule avec condition de Neumann $$\lambda_0^N(C) = \lambda_0.$$

Nous obtenons facilement une première minoration du spectre d'une surface modelée sur un graphe :

\begin{prop}\label{prop:MinGene}
Pour toute surface $M$ modelée sur le graphe $G$ à partir de la cellule $C$, on a 
$$\lambda_0(M)\geq \lambda_0.$$
\end{prop}

\begin{lemm}\label{lemm:MinGene}
Pour toute partie $M_f$ de $M$ contenue dans un nombre fini de cellules, son bas du spectre avec condition de Dirichlet vérifie
$$\lambda_0^D(M_f)\geq \lambda_0^N(C)= \lambda_0.$$
\end{lemm}
\begin{proof}
Il suffit de modifier légèrement la démonstration du Lemme \ref{lemm:MonoNeum} pour obtenir le résultat suivant : si $M$ est la réunion de $M_1\cup\ldots\cup M_k$ dont les intérieurs sont disjoints, alors 
$$\lambda_0^N(M)\geq \min_k(\lambda_0^N(M_1),\ldots,\lambda_0^N(M_k)).$$
Soit $\widetilde M$ une réunion d'un nombre fini de copies de $C$ d'intérieurs dsijoints, telle que $M_f\subset \widetilde M$, d'après l'argument précédent
$$\lambda_0^N(\widetilde M)\geq \lambda_0.$$
Comme $\Cl H_0^1(M_f)\subset \Cl H^1(\widetilde M)$ (voir Section \ref{ssec:PreSpec} pour les notations), on a 
$$\lambda_0^D(M_f)\geq\lambda_0^N(\widetilde M)\geq\lambda_0.$$
\end{proof}
Comme, par définition, on a 
$$\lambda_0(M) = \inf_{M_f\subset M} \lambda_0^D(M_f)$$ où $M_f$ parcourt l'ensemble des parties compactes de $M$, la Proposition \ref{prop:MinGene} découle immédiatement du Lemme \ref{lemm:MinGene}.

On obtient en particulier :

\begin{coro}\label{coro:L0Pos}
$\lambda_0(M)>0$
dès que $\lambda_0(C)>0$.
\end{coro}

\begin{rema}\label{rem:Double}
Le bas du spectre de $C$ avec condition de Neumann est exactement le bas du spectre de la surface sans bord $C^2$, réunion de deux copies de $C$ le long de $\bd C$ (voir Fig. 4) : ce résultat général de théorie spectrale s'obtient en combinant le Lemme \ref{lemm:MonoNeum} et la Proposition \ref{prop:InvEigFunc} (ou le Corollaire \ref{coro:InvEigFunc} si $C$ n'a pas de fonction propre associée à $\lambda_0$).
\end{rema}
D'après le Théorème \ref{th:Don}, des cellules $C$ dont le bas du spectre est positif existent : il suffit qu'elles soient géométriquement finies de volume infini, avec $\lambda_0$ suffisamment petit. On peut alors utiliser les résultats de la Section \ref{sec:GFin} pour contrôler le bas du spectre de $C^2$ en fonction du volume de son c\oe ur convexe. 

\begin{figure}\label{fig:double}
\begin{center}
\includegraphics[width = 0.6\textwidth]{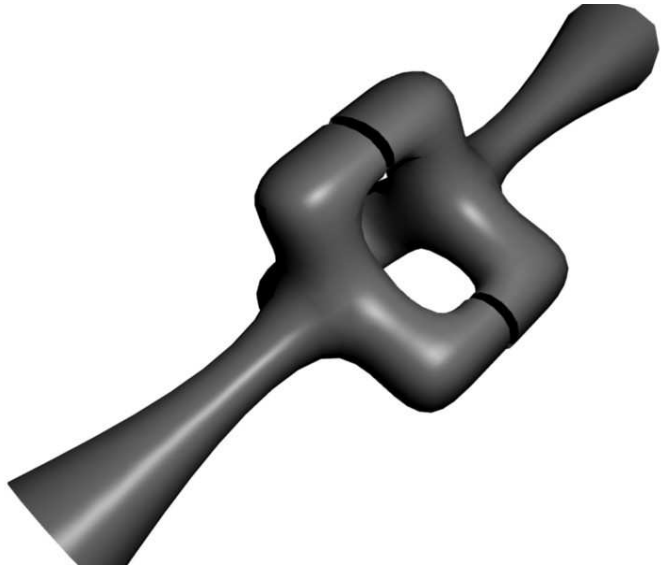}
\caption{Double d'une cellule à 3 composantes de bord}
\end{center}
\end{figure}

\subsection{Graphes moyennables et majoration de $\lambda_0$}\label{ssec:GInfMoy}

\begin{defi}
Soit $G$ un graphe infini, on appelle \emph{constante de Cheeger} de $G$ la constante
$$\g h(G) = \inf_{G_f\subset G} \frac{\#\bd G_f}{\#G_f},$$
où $G_f$ parcourt l'ensemble des parties finies de $G$ et $\bd G_f$ est l'ensemble des points de $G_f$ reliés à un point de $G\bs G_f$.
\end{defi}

On dira qu'un graphe est \emph{moyennable} si et seulement si $\g h(G) = 0$. Cette terminologie vient du résultat classique suivant :
\begin{theo}[F\o lner]
Soit $(\Gamma,S)$ un groupe de type fini muni d'un système de générateurs, $\Gamma$ est moyennable, c'est-à-dire admet une mesure de probabilité simplement additive invariante à gauche, si et seulement si le graphe de Cayley de $\Gamma$ relativement au système de générateurs $S$ est moyennable au sens défini précédemment pour les graphes.
\end{theo}

On peut par exemple trouver la démonstration de ce théorème dû à F\o lner dans \cite{Br}.

\pgh

Sur un graphe $G$, nous appellerons $\Delta_G$ le Laplacien combinatoire associé de la même façon qu'à la Section \ref{ssec:PincGeod} à la forme quadratique
$$q = \sum_{i\sim j}(f(i)-f(j))^2$$ et à la mesure
$$\mu = \sum_{i\in G}\delta(i).$$ Il s'agit alors de l'opérateur linéaire qui a toute fonction $f$ sur le graphe associe la fonction $\Delta_Gf$ définie par
$$\Delta_G f(i) = \sum_{i\sim j} (f(i) - f(j)).$$
C'est un opérateur autoadjoint positif, dont la première valeur propre $\mu_0(G)$ est l'infimum des quotients de Rayleigh combinatoires : 
$$\mu_0(G) = \inf_{f}\frac{\sum_{i\sim j}(f(i)-f(j))^2}{\sum_i f(i)^2}.$$

La moyennabilité du graphe est reliée à la première valeur propre de cet opérateur :
\begin{theo}[Inégalités de Cheeger combinatoires]\label{th:CheegComb}
$$\frac{1}{2v}\g h(G)^2\leq \mu_0(G)\leq \g h(G).$$
\end{theo}

En particulier, le graphe est moyennable si et seulement si la première valeur propre de son Laplacien combinatoire est nulle. On peut consulter \cite{CV} p31 pour une preuve de ce résultat que nous utiliserons ultérieurement.

L'objectif de cette section est de démontrer le théorème suivant :

\begin{theo}\label{th:GInfMoy}
Soit $M$ une surface hyperbolique modelée sur un graphe $G$ à valence constante à partir de la cellule $C$. On a
$$\lambda_0(M) \leq \lambda_0 + A_2 \g h(G),$$
où 
$$A_2 = (v-1)(\frac{1}{m(l)^2} + \lambda_0).$$
\end{theo}

On rappelle que $m(l)$ a été défini à la Section \ref{ssec:PreSurf} comme la largeur du voisinage collier d'une géodésique de longueur $l$, qui vaut
$$m(l) = {\rm argsinh}(\frac{1}{\sh(l/2)}).$$

En particulier, ceci prouve la première moitié du Théorème \ref{th:GInf1} :

\begin{coro}
Sous les hypothèses précédentes, si $G$ est moyennable, alors $$\lambda_0(M) = \lambda_0.$$
\end{coro}

\begin{proof}
Soit $G$ un graphe infini moyennable de valence constante $v$, $C$ une surface hyperbolique à bord compact, invariante par une isométrie $J$ d'ordre $v$. Soit $M$ la variété modelée sur $G$ à partir de $C$.

Soit $(G_n)_{n\in\Bb N}$ une famille de parties finies (que nous supposerons connexes) de $G$ telle que 
$$\lim_{n\ra \infty}\frac{\#\bd G_n}{\#G_n}=\mu_0(G),$$
une telle famille existe par définition d'un graphe moyennable. 
On note $M_n$ la réunion des cellules correspondant à la partie $G_n$, et $\bd^G M_n$ les cellules de $M\bs M_n$ reliées par au moins une composante de bord à une cellule de $M_n$. Il y a donc au plus $(v-1)\#\bd G_n$ cellules dans $\bd^G M_n$. On note $M_n^+ = M_n\cup \bd^G M_n$.

\pgh
Pour $\epsilon\in[0,1]$, on va construire une famille de fonctions $\tilde f_{\epsilon,n}$ à support compact dans $M_n^+$ telles que 
$$\frac{||\nabla \tilde f_{\epsilon,n}||^2}{||\tilde f_{\epsilon,n}||^2}\leq \lambda_0+\epsilon+(v-1)(\frac{1}{m(l)^2} + \lambda_0+\epsilon)\frac{\#\bd G_n}{\#G_n},$$
ce qui montrera le Théorème \ref{th:GInfMoy}. L'idée de cette démonstration est de reproduire sur toutes les cellules de $M_n$ une fonction dont le quotient de Raileygh est presque $\lambda_0$ et de la rendre continue à support compact sur $M_n^+$. Cette dernière opération ajoutera à l'énergie de $\tilde f_{\epsilon,n}$ de l'ordre de $\#\bd G_n$.

\pgh

Soit $f_\epsilon$ une fonction $\Cl C^\infty$ sur $C$, vérifiant les conditions de Neumann sur $\bd C$ et vérifiant
$$\int_C|\nabla f_\epsilon|^2 = (\lambda_0+\epsilon)\int_C f_\epsilon^2.$$
D'après la Section \ref{ssec:PreSpec}, on peut supposer que $f_\epsilon$ est invariante par $J$, et on supposera que
$$\int_C f_\epsilon^2 = 1.$$
On étend $f_\epsilon$ en $\tilde{f}_{\epsilon,n}$ définie sur $M_n$ par
$$\tilde{f}_\epsilon(x) = f_\epsilon\circ \phi_p(x)\mbox{ si } x\in C_{p}.$$
On rappelle que $\phi_p$ est l'isométrie qui identifie $C_p$ à $C$ ; désormais nous sous-entendrons cette identification en écrivant
$$\tilde{f}_\epsilon|_{C_p}(x) = f_\epsilon(x).$$
Comme $f$ est invariante par $J$, $f$ est identique sur chaque composante de bord. $\tilde{f}_\epsilon$ est donc continue sur $M_\epsilon$ et $\Cl C^2$ par morceaux, donc 
$\tilde{f}_\epsilon\in\Cl H^1(M_n).$ Il s'agit maintenant de lui donner un support compact dans $M_n^+$.

Soit $\alpha$ une géodésique fermée de longueur $l$ de $M_n\cap\bd^G M_n$, et $C_\alpha$ la cellule de $M\bs M_n$ à laquelle elle appartient. D'après la Proposition \ref{prop:LemCol}  le collier
$$\{p\in X : d(p,x)\leq m(l)\}$$
est un voisinage tubulaire de $\alpha$ plongé dans $X$, avec $$m(l) = {\rm argsinh}(\frac{1}{\sh(l/2)}).$$
Il admet des coordonnées $(r,\theta)$, où $|r|\leq m(l)$ et $\theta\in\Bb S^1$, telles que la métrique s'écrit 
$$ds^2 = dr^2+(\frac{l}{2\pi})^2\ch^2rd\theta^2.$$

On suppose que $r\leq0$ correspond à la partie du collier située dans $C_\alpha\subset\bd^G M_n\bs M_n$. Pour $r\leq0$, posons $$\psi_\alpha(r,\theta) = \frac{1}{m(l)}(m(l)-r).$$
Pour toute géodésique $\beta\subset C_\alpha\cap M_n$, on définit $\psi$ de la même façon sur le tube de $C_\alpha$ qui entoure $\beta$, et on l'étend sur $C_\alpha$ par $\psi_\alpha = 0$  à l'extérieur de ces voisinages tubulaires. $\psi_\alpha$ est continue, $\Cl C^2$ par morceaux, toujours inférieure à $1$, vaut $1$ sur les géodésiques de $C_\alpha\cap M_n$ et pour tous $r\in[0,m(l)],$
$$|\nabla \psi_\alpha (r,\theta)| = \frac{1}{m(l)}<\infty.$$

On a alors \beq \label{eq:EneCutOf} \int_{C_\alpha}|\nabla (\psi_\alpha f_\epsilon)|^2 = \int_C|\frac{f_\epsilon}{m(l)^2} + \psi\nabla f_\epsilon|^2\leq \frac{1}{m(l)^2}+\lambda_0+\epsilon\eeq 
car $||f_\epsilon||^2_{C_\alpha} = 1$.

Pour toute cellule $C_\alpha\subset \bd^G M_n$, on pose
$$\tilde{f}_{\epsilon,n} = \psi_\alpha f_\epsilon.$$

Nous obtenons finalement une fonction  $\tilde{f}_{\epsilon,n}$ continue et $\Cl C^1$ par morceaux à support compact dans l'intérieur de $M_n^+$, donc $\tilde{f}_{\epsilon,n}\in\Cl H^1_0(M_n^+)\subset \Cl H^1(M)$, qui vérifie : 

$$\int_{M_n^+}|\nabla \tilde{f}_\epsilon|^2 = \sum_{C_i\subset M_n^+}\int_{C_i}|\nabla \tilde{f}_\epsilon|^2\leq (\lambda_0+\epsilon)\sum_{C_i\subset M_n}\int_{C_i}\tilde{f}_\epsilon^2 + \sum_{C_i\subset \bd^G M_n}\int_C |\nabla (\psi f_\epsilon)|^2$$
d'où 
$$\int_{M_n}|\nabla \tilde{f}_\epsilon|^2\leq \# G_n(\lambda_0+\epsilon) + (v-1)\#(\bd G_n)(\frac{1}{m(l)^2}+\lambda_0+\epsilon)$$
d'après (\ref{eq:EneCutOf}).

On obtient finalement
$$\frac{\int_{M_n}|\nabla \tilde{f}_\epsilon|^2}{\int_{M_n}\tilde{f}_\epsilon^2}\leq \frac{\#G_n(\lambda_0+\epsilon)||f_\epsilon||^2+(v-1)\#\bd G_n K'||f_\epsilon||^2}{\#G_n||f_\epsilon||^2}\leq \lambda_0+\epsilon + (v-1)(\frac{1}{m(l)^2} + \lambda_0+\epsilon)\frac{\#\bd G_n}{\#G_n},$$
ce qui conclut notre preuve.
\end{proof}

\begin{rema} Cette partie de notre démonstration n'utilise pas l'existence de la fonction propre $\psi_0$, elle est donc valable y compris lorsque la cellule n'a pas un trou spectral positif. Elle est analogue à celle de \cite{Br},§2 et a été utilisée depuis dans de nombreux articles traitant de moyennabilité.
\end{rema}

\begin{rema}
Dans le cas où $G$ est le graphe de Cayley d'un groupe abélien (par exemple $\Bb Z, \Bb Z^2...$), donc moyennable, on retrouve ce qui est appelé communément une \emph{surface périodique}. Le résultat $\lambda_0(M) = \lambda_0(C)$ est alors un corollaire immédiat de la \emph{théorie de Floquet}. Cependant, d'une part notre démonstration (d'un résultat certes beaucoup plus faible) est plus élémentaire que la construction de la théorie de Floquet, d'autre part elle s'applique ici à une classe de surfaces beaucoup plus large, que nous appelons également \emph{périodiques} au sens de la Définition \ref{def:SInfPer}. Nous reparlerons plus en détails de la situation où $G$ est un graphe de Cayley au Paragraphe \ref{ssec:GInfRev}.
\end{rema}

\begin{coro}
Si $G$ est moyennable, le bas du spectre de la surface $M$ est le même que celui de $C^2$, où $C^2$ est la surface hyperbolique non compacte complète sans bord, \emph{double} de $C$ (voir Fig. \ref{fig:double}).
\end{coro}

\begin{proof}
Il suffit de se rappeler que d'après la Remarque \ref{rem:Double}, le bas du spectre de $C$ avec condition de Neumann est le même que celui de son double $C^2$. 
\end{proof}
Nous avons donc toute une famille de surfaces hyperboliques dont le bas du spectre est égal à celui de la surface $C^2$. Comme cas particuliers des graphes moyennables, citons les graphes finis et les graphes à croissance polynomiale, parmi lesquels les graphes abéliens cités précédemment. Notons qu'il existe des graphes à croissance exponentielle qui restent moyennables.

Le Théorème \ref{th:GInfMoy} nous donne la première implication du Théorème \ref{th:GInf1}. Nous allons maintenant nous intéresser à l'autre.

\subsection{Graphes non moyennables et minoration de $\lambda_0$}

Soit $G$ un graphe non moyennable de valence constante $v$, de constante de Cheeger $\g h_G>0$, $C$ une surface hyperbolique à bord compact invariant par une isométrie $J$ d'ordre $v$ et $M$ la surface modelée selon $G$ à partir de $C$. On note $\mu_0(G)$ le bas du spectre du Laplacien combinatoire sur $G$ (voir Théorème \ref{th:CheegComb}). Cette section sera consacrée à la démonstration du théorème suivant :

\begin{theo}\label{th:GInfNMoy}
Avec les notations précédentes, on a
$$\lambda_0(M)\geq\lambda_0^N(C) +A_1\mu_0,$$
où $A_1$ dépend de caractéristiques spectrales de $C$ et de la longueur de $\bd C$.
\end{theo}

Ce résultat nous donne la deuxième inégalité du Théorème \ref{th:GInf2} et prouve donc à l'aide du Théorème \ref{th:CheegComb} la deuxième implication du Théorème \ref{th:GInf1}.

\pgh

Par hypothèse, $C$ possède une unique fonction propre positive $\psi_0$ associée à $\lambda_0 = \lambda_0(C)$, que nous supposerons de norme $1$ :
$$\int_C\psi_0^2 = 1,$$
invariante par $J$, et a un trou spectral $\eta = \lambda_1-\lambda_0>0$. Notons
$$\tilde \lambda_0 = \lambda_0(M)\mbox{, et } \delta = \tilde\lambda_0-\lambda_0.$$
D'après la Proposition \ref{prop:MinGene}, $\delta\geq 0$.

\begin{proof}
Soit $(f_\epsilon)_{0<\epsilon<1}$ une famille de fonctions à support compact dans $M$ telle que pour tout $\epsilon>0$, on ait
$$\frac{\int_M|\nabla f_\epsilon|^2}{\int_M f_\epsilon^2}\leq \tilde\lambda_0 + \epsilon.$$

Nous allons discrétiser ces fonctions en les projetant, cellule par cellule, sur la fonction propre $\psi_0$. La composante perpendiculaire à $\psi_0$, nécessairement non nulle pour obtenir un support compact, empêchera d'avoir $\tilde\lambda_0 = \lambda_0$. Pour chaque cellule $C_i$ de $M$, posons
$$a_i^2(\epsilon) = ||f_\epsilon||^2_{C_i} = \int_{C_i}f_\epsilon^2,$$
$$b_i(\epsilon) = \langle f_\epsilon, \psi_0\rangle_{C_i} = \int_{C_i}f_\epsilon\psi_0,$$
et
$$c_i^2(\epsilon) = ||f_\epsilon-b_i\psi_0||^2_{C_i} = \int_{C_i}(f_\epsilon-b_i\psi_0)^2.$$
On a alors $a_i^2 = b_i^2+c_i^2$. Nous noterons désormais $$g_i = f_\epsilon-b_i\psi_0,$$ où $g_i$ est la composante de $f_i$ orthogonale à $\psi_0$ pour le produit scalaire $\Cl L^2$.

Sur chaque cellule $C_i$,
$$\int_{C_i} \psi_0 g_i = \int_{C_i} \psi_0 f_\epsilon - b_i \int_{C_i} \psi_0^2 = 0$$ et $$\int_{C_i} \nabla \psi_0.\nabla g_i = \int_{C_i}\Delta \psi_0 g_i = \lambda_0\int_{C_i} \psi_0 g_i = 0.$$
On a donc
$$\frac{\int_M |\nabla f_\epsilon|^2}{\int_M f_\epsilon^2} = \frac{\sum_i ||\nabla f_\epsilon||^2_{C_i}}{\sum_i || f_\epsilon||^2_{C_i}} = \frac{\sum_i b_i^2||\nabla \psi_0||^2_{C_i}+||\nabla g_i||^2_{C_i}}{\sum_i a_i^2} = \frac{\sum_i \lambda_0b_i^2+||\nabla g_i||^2_{C_i}}{\sum_i a_i^2}.$$
Comme $g_i$ est orthogonale à $\psi_0$, $$||\nabla g_i||^2_{C_i}\geq \lambda_1 ||g_i||^2_{C_i} = \lambda_1 c_i^2.$$ On obtient alors 
$$ ||\nabla g_i||^2_{C_i} = \frac{\lambda_0}{\lambda_1}||\nabla g_i||^2_{C_i}+\frac{\eta}{\lambda_1}||\nabla g_i||^2_{C_i}\geq \lambda_0c_i^2+\frac{\eta}{\lambda_1}||\nabla g_i||^2_{C_i},
$$
d'où
$$\tilde \lambda_0+\epsilon\geq \frac{\int_M |\nabla f_\epsilon|^2}{\int_M f_\epsilon^2}\geq \frac{\sum_i \lambda_0b_i^2+\lambda_0c_i^2+\frac{\eta}{\lambda_1}||\nabla g_i||^2_{C_i}}{\sum_i a_i^2} = \lambda_0 + \frac{\eta}{\lambda_1}\frac{\sum_i||\nabla g_i||^2_{C_i}}{\sum_i a_i^2}.$$
On a donc en particulier
\beq \label{eq:MinNMoy} \delta +\epsilon = \tilde\lambda_0-\lambda_0\geq +\epsilon\geq \frac{\eta}{\lambda_1}\frac{\sum_i||\nabla g_i||^2_{C_i}}{\sum_i a_i^2}.\eeq

Nous allons donc nous intéresser au terme $\frac{\sum_i||\nabla g_i||^2_{C_i}}{\sum_i a_i^2}$ et montrer le lemme :

\begin{lemm}\label{lemm:AnaComb}
Il existe une constante $A$ ne dépendant que de propriétés spectrales de la cellule $C$ et de la longueur de $\bd C$ telle que
$$\sum_i||\nabla g_i||_{C_i}^2\geq A\sum_{i\sim j}(b_i-b_j)^2.$$
\end{lemm}

\begin{proof}

Soit $i\sim j$ et $\alpha_{ij}$ la géodésique commune à deux cellules $C_i$ et $C_j$, elle est de longueur $l$ comme toutes les géodésiques du bord des cellules (toutes isométriques). D'après la Section \ref{ssec:PreSurf}, soit $m=m(l)$, il existe dans $C_i$ et dans $C_j$ un voisinage tubulaire $T_{ij}$ de $\alpha_{ij}$ sur lequel la métrique hyperbolique s'écrit $$ds^2 = dr^2 + (\frac{l}{2\pi})^2\ch^2r d\theta^2$$ pour $-m\leq r\leq m$, où par convention $r\geq 0$ si et seulement si $(r,\theta)\in C_i$. On note $T_{ij}^+$ la partie de $T_{ij}$ située dans $C_i$, et $T_{ij}^-$ l'autre. 

Puisque nous voulons une minoration, nous pouvons nous limiter à l'étude de ce qu'il se passe sur ces tubes :

\beq \label{eq:MinLoc} \sum_i||\nabla g_i||^2_{C_i}\geq \frac{1}{2}\sum_i(||\nabla g_i||^2_{C_i}+\lambda_1||g_i||^2_{C_i})\geq \sum_{i\sim j}(||\nabla g_i||^2_{T_{ij}^+}+||\nabla g_j||^2_{T_{ij}^-}+\lambda_1(||g_i||^2_{T_{ij}^+}+||g_j||^2_{T_{ij}^-})).\eeq

Nous poursuivons alors notre minoration par le lemme suivant :

\begin{lemm}\label{lemm:AnaComb2}
Il existe une constante $A$ ne dépendant que de $l$ et de propriétés spectrales de la cellule $C$ telle que pour tous $(i,j)\in G^2$ avec $i\sim j$,  
$$||\nabla g_i||^2_{T_{ij^+}}+||\nabla g_j||^2_{T_{ij}^-}+\lambda_1(||g_i||^2_{T_{ij}^+}+||g_j||^2_{T_{ij^-}})\geq A(b_i-b_j)^2.$$
\end{lemm}

\begin{proof}
On note $f(r,\theta)$ l'expression de $f$ en coordonnées de Fermi sur le tube $T_{ij}$. 

Pour nous ramener à un problème ne dépendant que de $r$, pour chaque fonction $f$ définie sur $T_{ij}$, notons 
$$F(r) = \frac{1}{2\pi}\int_0^{2\pi}f(r,\theta)d\theta.$$ Remarquons que d'après l'inégalité de Cauchy-Schwarz, 

$$\left(\int_0^{2\pi}f(\theta)d\theta\right)^2\leq(2\pi)^2\int_0^{2\pi}f^2(\theta)d\theta.$$
On a donc
$$\int_{T_{ij}}f^2 = \int_{-m}^{m}\int_0^{2\pi} f(r,\theta)^2\frac{l}{2\pi}\ch rd\theta dr\geq \int_{-m}^{m}\frac{(\int_0^{2\pi}fd\theta)^2}{(2\pi)^2} \frac{l}{2\pi}\ch rdr = \frac{1}{2\pi}\int_{T_{ij}}F^2.$$
De plus,
$$ \int_{T_{ij}}|\nabla f|^2 = \int_{-m}^{m}\int_0^{2\pi}[(\frac{\bd f}{\bd r})^2+\frac{1}{(l/2\pi)^2\ch^2r}(\frac{\bd f}{\bd\theta})^2 ]\frac{l}{2\pi}\ch rd\theta dr\geq \int_{-m}^{m}\int_0^{2\pi}(\frac{\bd f}{\bd r})^2\frac{l}{2\pi} \ch rd\theta dr.$$
On obtient alors, de même que précédemment,
$$\int_{T_{ij}}|\nabla f|^2\geq \int_{-m}^{m}\frac{(\int_0^{2\pi}\frac{\bd f}{\bd r}d\theta)^2}{(2\pi)^2} \frac{l}{2\pi}\ch rdr = \frac{1}{2\pi}\int_{T_{ij}}|\nabla F|^2.$$

Pour montrer le Lemme \ref{lemm:AnaComb2}, il suffit donc de montrer le même résultat pour les $G_i$ (moyennes cylindriques  de $g_i$), c'est-à-dire montrer qu'il existe $A'= 2\pi A$ ne dépendant que de $C$ tel que
$$||G'_i||^2_{T_{i,j}^+}+||G'_j||^2_{T_{i,j}^-}+\lambda_1(||G_i||^2_{T_{i,j}^+}+||G_j||^2_{T_{i,j}^+})\geq A''(b_i-b_j)^2.$$
Comme $\psi_0$ est invariante par $J$, elle est identique sur toutes les composantes de bord. Par continuité de $f_\epsilon$, on donc
$$g_j(0,\theta)+b_j \psi_0(0,\theta) = g_i(0,\theta)+b_i \psi_0(0,\theta),$$
d'où $$|G_j(0)-G_i(0)| = |b_j-b_i| \Psi_0(0).$$
Supposons, quitte à inverser $i$ et $j$, que $|G_i(0)|\geq \frac{1}{2}|b_j-b_i| \Psi_0(0),$ et notons
$$R = \inf\{r\in]0,m] : |G_i(r)|\leq \frac{1}{4}|b_j-b_i| \Psi_0(0)\}.$$

\emph{$1^{er}$ cas : R = m}\\
On a alors $$\forall r\in[0,m], |G_i(r)|\geq \frac{1}{4}|b_j-b_i| \Psi_0(0),$$
donc $$\lambda_1||G_i||^2_{T_{i,j}^+} = \lambda_1\int_{T_{ij^+}}| G_i|^2\geq \lambda_1\frac{1}{4}l\sh(m(l))|\Psi_0(0)|^2(b_j-b_i)^2=A''|\Psi_0(0)|^2(b_j-b_i)^2,$$
avec $$A'' = \lambda_1\frac{1}{4}l\sh(m(l)).$$

\emph{$2^{nd}$ cas : R<m}\\
Posons 
$$P= G_i(0) \geq \frac{1}{2}|b_j-b_i|\Psi_0(0)$$ et 
$$Q = G_i(R) = \frac{1}{4}|b_j-b_i|\Psi_0(0).$$ 
Soit $F_0$ la fonction sur le collier $[0,R]\cx[0,l]$ muni de la métrique de Fermi, ne dépendant que de $R$, valant $P$ en $0$, $Q$ en $R$ et harmonique pour le Laplacien hyperbolique : on vérifie à l'aide de (\ref{eq:LapFermi}) que l'on a

$$F_0(r) = P-\frac{P-Q}{U(R)}U(r)$$ avec $$U(r) = \arcsin(\Th r).$$

La fonction $F_0$ minimise alors l'énergie de Dirichlet (associée à la métrique de Fermi) parmi toutes les fonctions $H$ sur $\Bb S^1\cx[0,R]$ vérifiant $H(\theta,0) = P$ et $H(\theta,R) = Q$. On a donc nécessairement
$$||G'_i||^2_{T_{i,j}^+}\geq ||G'_i||^2_{\Bb S^1\cx[0,R]}\geq||F'_0||^2_{\Bb S^1\cx[0,R]} = \frac{(P-Q)^2}{U(R)^2}\int_0^R (U'(r))^2\frac{l}{2\pi}\ch r dr.$$

Or, $$U'(r) = \frac{1}{\ch^2r\sqrt{1-th^2r}},$$
donc $$U'(r)^2 = \frac{1}{\ch^4 r(1-th^2r)} = \frac{1}{\ch^2r},$$
d'où $$\int_0^R (U'(r))^2\frac{l}{2\pi}\ch r dr = \int_0^R \frac{l}{2\pi}\frac{1}{\ch r} dr = \frac{l}{\pi}\arctan(e^R)-\pi/2.$$
On a alors 
$$||G'_i||^2_{T_{i,j}^+}\geq (P-Q)^22\frac{2\arctan(e^R)-\pi/2}{\pi U(R)^2}\geq \frac{(b_i-b_j)^2|\Psi_0(0)|^2}{8}\frac{2\arctan(e^R)-\pi/2}{\pi U(R)^2},$$
car $$|P-Q|\geq \frac{b_i-b_j}{4}|\Psi_0(0)|.$$ 

\pgh

Un développement limité nous montre que 
$$\frac{2\arctan(e^R)-\pi/2}{U(R)^2}\sim\frac{1}{R}$$
lorsque $R$ tend vers $0$, il existe donc une constante $A'''$, qui ne dépend que de $l$, telle que pour tout $R\in]0,m(l)]$, on ait
$$\frac{2\arctan(e^R)-\pi/2}{U(R)^2}>\frac{8\pi}{l}A'''.$$

On a donc $$||G'_i||^2_{T_{i,j}^+}>A'''|\Psi_0(0)|^2(b_i-b_j)^2.$$

Finalement, posons $$A' = \max(A'', A''')|\Psi_0(0)|^2$$ qui ne dépend que de $\lambda_1$, de $\Psi_0(0)$ et de $l$, on a
$$||G'_i||^2_{T_{i,j}^+}+||G'_j||^2_{T_{i,j}^-}+\lambda_1(||G_i||^2_{T_{i,j}^+}+||G_j||^2_{T_{i,j}^+})\geq A'(b_i-b_j)^2,$$
ce qui conclut la démonstration de notre lemme.
\end{proof}

La minoration (\ref{eq:MinLoc}) devient alors

$$\sum_i||\nabla g_i||^2\geq A\sum_{i\sim j}(b_i-b_j)^2,$$
avec $A = \frac{A'}{2\pi}$ qui dépend de $\lambda_1$, de $\Psi_0(0)$ et de $l$, ce qui conclut la démonstration du Lemme \ref{lemm:AnaComb}.
\end{proof}

On a d'après l'inégalité (\ref{eq:MinNMoy}),
$$\frac{\eta}{\lambda_1}\sum_i||\nabla g_i||^2_{C_i}\leq (\delta+\epsilon)||f_\epsilon||^2 = (\delta+\epsilon)\sum_i a_i^2,$$
donc
$$\eta\sum_i||g_i||^2_{C_i} = \eta\sum_i c_i^2\leq (\delta+\epsilon)\sum_i a_i^2.$$
On a alors
$$\sum_i a_i^2 = \sum_i b_i^2+\sum_i c_i^2\leq \sum_i b_i^2+\epsilon/\eta\sum_i a_i^2,$$
d'où $$\sum_i a_i^2\leq \frac{1}{1-\frac{\delta+\epsilon}{\eta}}\sum_i b_i^2.$$
L'inégalité (\ref{eq:MinNMoy}) devient
$$\delta + \epsilon\geq\frac{\eta}{\lambda_1}\frac{\sum_i||\nabla g_i||^2_{C_i}}{\sum_i a_i^2}\geq (1-\frac{\delta+\epsilon}{\eta})\frac{\eta A}{\lambda_1}\frac{\sum_{i\sim j}(b_i-b_j)^2}{\sum_i b_i^2}.$$
On a vu que $$\mu_0 = \inf \frac{\sum_{i\sim j}(\alpha_i-\alpha_j)^2}{\sum_i \alpha_i^2},$$ où $(\alpha_i)$ parcourt l'ensemble des familles positives à support compact dans $G$ est le bas du spectre du Laplacien combinatoire de $G$. D'après le Théorème \ref{th:CheegComb}, lorsque G n'est pas moyennable, $\mu_0>0$. 

On obtient donc en faisant tendre $\epsilon$ vers $0$ :
$$\delta\geq \frac{\eta-\delta}{\lambda_1}A\mu_0,$$
soit finalement
\beq \label{eq:MinDelta} \delta\geq \frac{\eta}{1+1/\lambda_1}A\mu_0.\eeq
Ceci conclut la démonstration du Théorème \ref{th:GInfNMoy}, avec 
$$A_1 = \frac{\eta}{1+1/\lambda_1}A,$$
qui dépend de $\lambda_1$, $\eta$, $\Psi_0(0)$ et $l$.
\end{proof}

\begin{rema}
Cette deuxième partie de la démonstration utilise de façon essentielle l'hypothèse $\lambda_1>\lambda_0$.
\end{rema}

\section{Quelques généralisations}

Nous remarquons que dans nos démonstrations, l'invariance de la cellule par une isométrie d'ordre fini n'est utilisée que pour obtenir l'invariance de la première fonction propre et pouvoir recoller nos fonctions d'une cellule à l'autre. Nous présentons ici deux situations où l'isométrie cyclique de la cellule ne sera pas nécessaire, ce qui va nous permettre d'utiliser les mêmes méthodes avec des surfaces géométriquement infinies non périodiques. Il s'agit du cas où le volume des cellules est fini et uniformément borné, et du cas où notre cellule est en fait un domaine fondamental pour l'action d'un groupe de revêtement.

\subsection{Surfaces à découpage borné}\label{ssec:GInfBorn}

Si la cellule est de volume fini, la première fonction propre est constante, $\lambda_0 = 0$ et son trou spectral est positif. Nous obtenons alors un cadre très simple où utiliser les méthodes de la Section \ref{sec:GInf} :

\begin{defi}\label{def:DecBorn}
Soit $M$ une surface hyperbolique, on dira que $M$ admet un \emph{découpage borné} s'il existe des constantes $\eta,v,k,K>0$ et une famille $(M_i)_i$ de sous-surfaces de $M$, d'intérieurs disjoints, à bords géodésiques, telles que $\forall i$, le nombre de composantes connexe de $\bd M_i$ est borné par $v$, le trou spectral $\lambda_1(C_i)>\eta$,
$$k<\Vol(M_i)<K$$ 
et pour toute composante de bord $\alpha\subset\bd M_i$, $$k<\ell(\alpha)<K.$$
\end{defi}

Soit $M$ une variété admettant un découpage borné. On considère alors le graphe $G = (V,E)$ suivant : $V$ est l'ensemble des composantes $M_i$ du découpage de $M$, et si $(i,j)\in V^2, i\neq j,(i,j)\in E$ si et seulement si $M_i\cap M_j=\alpha_{ij}\neq \vd$. Nous appellerons $G$ le graphe sous-jacent au découpage de $M$.

Soit $M$ une surface hyperbolique admettant un découpage borné de constantes $\eta, v,k,K$ dont le graphe sous-jacent est $G$. Les deux propositions qui suivent démontrerons le Théorème \ref{th:GInfBorn} :

\begin{prop}
Sous les hypothèses ci-dessus, 
$$\lambda_0(M) \leq \frac{K(v-1)}{m(K)^2}\g h(G).$$
En particulier, si $G$ est moyennable, $\lambda_0(M) = 0$.
\end{prop}

\begin{proof}
Il suffit de reprendre la preuve du Théorème \ref{th:GInfMoy}, avec $f_\epsilon = 1$ pour tout $\epsilon>0$. Nous utilisons les notations de cette démonstration. La fonction cut-off $\psi_\alpha$ créée sur chaque tube de $\bd^G M_n$ dépend de la longueur $\ell(\alpha)$ de la géodésique du tube : l'inégalité (\ref{eq:EneCutOf}) devient donc simplement
$$\int_{C_\alpha} |\nabla \psi_\alpha|^2\leq \int_{C_\alpha}\frac{1}{m(\ell(\alpha))^2}\leq \frac{K}{m(K)},$$
car $l\fa m(l)$ est décroissante. La conclusion de la proposition s'ensuit.
\end{proof}

\begin{rema} 
Pour cette partie de notre démonstration, seules les hypothèses $v$ uniformément majoré et $$\ell(\alpha)\leq K \mbox{ et } \Vol(C_i)\leq K$$
pour toute cellule $C_i$ et toute composante de bord $\alpha$ d'une cellule sont nécessaires ; aucune condition spectrale n'est requise.
\end{rema}

\begin{prop}\label{prop:MinBorn}
Si $M$ admet un découpage borné dont le graphe sous-jacent est $G$, alors il existe une constante $A$ qui ne dépend que des constantes $v,\eta,k$ et $K$ telle que 
$$\lambda_0(M)\geq A\mu_0(G).$$
\end{prop}

\begin{proof}
Nous allons survoler, étape par étape, la démonstration du Théorème \ref{th:GInfNMoy} afin de voir comment l'adapter. Notre seule tâche est d'exprimer les constantes de minoration de cette preuve uniquement en fonction des constantes $v,\eta,k$ et $K$. Pour obtenir une normalisation des fonctions propres compatibles entre les différentes cellules, nous notons désormais pour toute fonction $h$ sur une cellule $C_i$ du découpage de $M$ : 
$$||h||^2_{\Cl L^2(C_i)} = \frac{1}{\Vol(C_i)}\int_{C_i}h^2,$$
et nous considérons le produit scalaire associé normalisé par $1/\Vol(C_i)$. Sur chaque cellule, la première fonction propre du Laplacien avec condition de Neumann de norme $1$ pour ce produit scalaire est donc $\psi_0 = 1$. 

\pgh

Pour tout $\epsilon>0$, soit $f_\epsilon$ de classe $\Cl C^\infty$ à support compact dans $M$ telle que 
$$\frac{||\nabla f_\epsilon||^2_{\Cl L^2(M)}}{||f_\epsilon||^2_{\Cl L^2(M)}}\leq \lambda_0(M)+\epsilon.$$

Nous posons alors pour chaque cellule $C_i$ de $M$ :

$$a_i^2(\epsilon) = ||f_\epsilon||^2_{C_i} = \frac{1}{\Vol(C_i)}\int_{C_i}f_\epsilon^2,$$
$$b_i(\epsilon) = \langle f_\epsilon, 1\rangle_{C_i} = \frac{1}{\Vol(C_i)}\int_{C_i}f_\epsilon,$$
$$c_i^2(\epsilon) = ||f_\epsilon-b_i||^2_{C_i} = \frac{1}{\Vol(C_i)}\int_{C_i}(f_\epsilon-b_i)^2.$$
On a alors $a_i^2 = b_i^2+c_i^2$. Nous notons encore $$g_i = f_\epsilon-b_i,$$ $g_i$ est la composante de $f_i$ orthogonale à $\psi_0$ pour le produit scalaire $\Cl L^2$. On a alors

\beq \label{eq:MinNMoyBorn} \frac{k\sum_i||\nabla g_i||^2_{C_i}}{K\sum_i a_i^2}\leq\frac{\sum_i||\nabla g_i||^2_{C_i}\Vol(C_i)}{\sum_i a_i^2\Vol(C_i)} = \frac{\int_M |\nabla f_\epsilon|^2}{\int_M f_\epsilon^2} \leq\lambda_0(M)+\epsilon = \tilde\lambda_0+\epsilon.\eeq

Nous allons donc de nouveau montrer l'analogue du Lemme \ref{lemm:AnaComb} :

\begin{lemm}\label{lemm:AnaCombBorn}
Il existe une constante $A$ ne dépendant que des constantes $\eta,v,k,K$ telle que
$$\sum_i||\nabla g_i||_{C_i}^2\geq A\sum_{i\sim j}(b_i-b_j)^2.$$
\end{lemm}

\begin{proof}
De même qu'en (\ref{eq:MinNMoy}), soit $i\sim j$ et $\alpha_{ij}$ la géodésique commune à $C_i$ et $C_j$ ; elle est de longueur $l\leq K$ . D'après la Section \ref{ssec:PreSurf}, il existe dans $C_i$ et dans $C_j$ un voisinage tubulaire $T_{ij}$ de $\alpha_{ij}$ sur lequel la métrique hyperbolique s'écrit $$ds^2 = dr^2 + (\frac{l}{2\pi})^2\ch^2r d\theta^2$$ pour $-m\leq r\leq m$, avec $m=m(K)\leq m(l_{\alpha_{ij}}$, où par convention $r\geq 0$ si et seulement si $(r,\theta)\in C_i$. On note $T_{ij}^+$ la partie de $T_{ij}$ située dans $C_i$, $T_{ij}^-$ l'autre. 

On a alors de nouveau
\beq \label{eq:MinLocBorn} \sum_i||\nabla g_i||^2_{C_i}\geq\sum_{i\sim j}(||\nabla g_i||^2_{T_{ij}^+}+||\nabla g_j||^2_{T_{ij}^-}+\lambda_1(||g_i||^2_{T_{ij}^+}+||g_j||^2_{T_{ij}^-})).\eeq

Nous poursuivons toujours notre minoration par le lemme suivant :

\begin{lemm}
Il existe une constante $A'$ ne dépendant que des constantes $\eta,v,k,K$ telle que pour tous $i\sim j$,  
$$||\nabla g_i||^2_{T_{ij^+}}+||\nabla g_j||^2_{T_{ij}^-}+\lambda_1(||g_i||^2_{T_{ij}^+}+||g_j||^2_{T_{ij}^-})>A'(b_i-b_j)^2.$$
\end{lemm}

\begin{proof}
Soit $i\sim j$, et $\alpha_{ij}$ une géodésique commune. Sur le voisinage tubulaire de $\alpha_{ij}$ de largeur $a = m(K)$, plongé dans $C_i\cup C_j$ quels que soient $i\sim j$, il suffit de recopier mot pour mot la démonstration du Lemme \ref{lemm:AnaComb2}. Celle-ci reste valable pour toutes les cellules car d'après l'expression (\ref{eq:LapFermi}), 
$$U(r) = \arcsin(\Th r)$$ est harmonique sur un voisinage tubulaire en coordonnées de Fermi, quelle que soit la longueur de la géodésique fermée qu'il entoure.
\end{proof}

On obtient donc d'après (\ref{eq:MinLocBorn}) 
$$\sum_i||\nabla g_i||^2_{C_i}\geq A\sum_{i\sim j}(b_i-b_j)^2,$$
ce qui conclut la démonstration du Lemme \ref{lemm:AnaCombBorn}.
\end{proof}

$$\frac{k\sum_i||\nabla g_i||^2_{C_i}}{K\sum_i a_i^2}\leq\frac{||\nabla f_\epsilon||^2}{||f_\epsilon||^2}\leq\tilde\lambda_0+\epsilon$$
 devient 
$$\sum_i c_i^2\leq \frac{K(\tilde\lambda_0+\epsilon)}{k\eta}\sum_i a_i^2,$$
donc
$$\sum_i a_i^2= \sum_i b_i^2+\sum_i c_i^2\leq\sum_i b_i^2+\frac{K(\tilde\lambda_0+\epsilon)}{k\eta}\sum a_i^2,$$
d'où
$$\sum_ia_i^2\leq \frac{1}{1-\frac{K(\tilde\lambda_0+\epsilon)}{k\eta}}\sum_ib_i^2.$$
On a donc d'après le Lemme \ref{lemm:AnaCombBorn}
$$\tilde\lambda_0 +\epsilon\geq \frac{\sum_i||\nabla g_i||^2_{C_i}}{\sum_i a_i^2}\geq (1-\frac{K(\tilde\lambda_0+\epsilon)}{k\eta})A\frac{\sum{i\sim j}(b_i-b_j)^2}{\sum_ib_i^2}\geq(1-\frac{K(\tilde\lambda_0+\epsilon)}{k\eta})A\mu_0(G).$$
Lorsque $\epsilon\ra 0$, on obtient de même qu'en (\ref{eq:MinDelta})
$$\lambda_0(M)\geq \frac{\eta}{1+K/(k\eta)}\mu_0(G),$$
ce qui conclut la preuve de la Proposition \ref{prop:MinBorn}.

\end{proof}

\subsection{Revêtements riemanniens et graphes de Cayley}\label{ssec:GInfRev}

\paragraph{Bas du spectre d'une surface modelée sur un graphe de Cayley}

Soit $G$ le graphe de Cayley d'un groupe de type fini $\Gamma$ associé aux générateurs $(g_1,...,g_n)$ qu'on suppose non triviaux. C'est un graphe à valence constante $v = 2n$, on peut donc à partir d'une cellule à $v$ composantes de bord invariante par une isométrie cyclique d'ordre $v$ considérer la surface $M$ modelée sur $G$ à partir de $C$. Le groupe $\Gamma$ agit canoniquement sur $G$ ; cette action se transporte naturellement sur $M$ : un élément $\gamma\in\Gamma$ envoie un point $x\in C_i$ sur le point 
$$\gamma(x) = \phi_{\gamma(i)}\circ(\phi_i)^{-1}(x),$$
en notant toujours $\phi_i$ l'isométrie qui identifie $C$ à $C_i$. Le quotient $M_1 = M/\Gamma$ est une surface sans bord, isométrique à la cellule $C$ dont les composantes de bord ont été identifiés deux à deux (Fig. 5): pour tout $1\leq p\leq n$, $$\alpha_{g_p} = C_i\cap C_{g_p(i)}$$ est identifié avec $$\alpha_{g_p^{-1}} = C_i\cap C_{g_p^{-1}(i)}.$$

\begin{figure}
\begin{center}
\includegraphics[width = 0.6\textwidth]{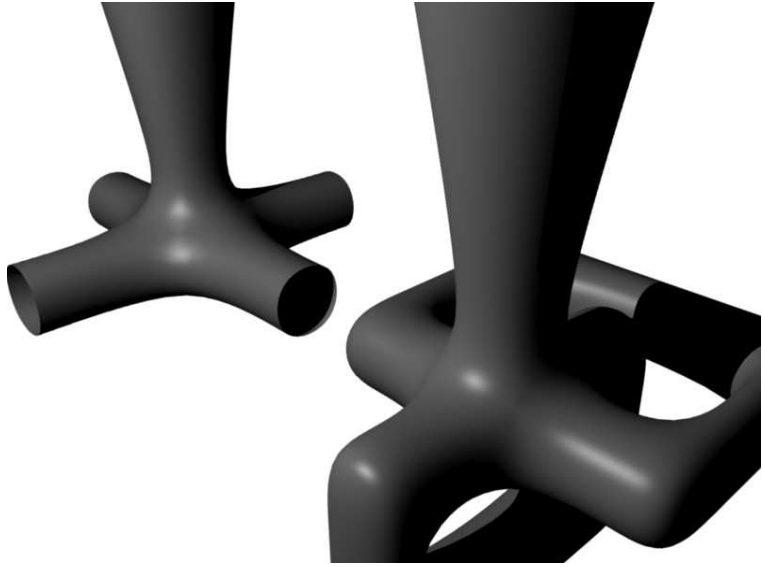}
\caption{Surface quotient $M/\Gamma$, à l'arrière plan la cellule $C$.}
\end{center}
\end{figure}

Le revêtement riemannien $M\ra M_1$ est galoisien de groupe $\Gamma$, et pour les bas des spectres de $M_1$ et de $C$ on a :

\begin{lemm}\label{lemm:SpecDomFond1}
Avec les notations précédentes, on a $$\lambda_0(M_1) = \lambda_0^N(C).$$
\end{lemm}

\begin{proof}
Par définition,
$$\lambda_0^N(C) = \inf\{\frac{||\nabla f||^2}{||f||^2}\}$$
où $f$ parcours les fonctions $\Cl C^\infty$ à support compact dans $C$, et 
$$\lambda_0(M_1) = \inf\{\frac{||\nabla f||^2}{||f||^2}\}$$
où $f$ parcours les fonctions $\Cl C^\infty$ à support compact dans $M_1$. Toute fonction $f\in\Cl C_0^\infty (M_1)$  se relève en une fonction $\tilde f\in\Cl C^\infty_0(C)$, et leurs quotient de Raileygh sont identiques. On a donc 
$$\lambda_0^N(C)\leq \lambda_0(M_1).$$

De plus, d'après le Théorème \ref{th:CarSul},
$$\lambda_0^N(C) = \sup\{\lambda\in\Bb R : \exists f\in\Cl C^\infty(C), f>0 : \Delta f = \lambda f\}$$
et 
$$\lambda_0(M_1) = \sup\{\lambda\in\Bb R : \exists f\in\Cl C^\infty(M_1), f>0 : \Delta f = \lambda f\}.$$
Toute fonction positive $\lambda$-harmonique sur $M_1$ se relève en une fonction positive $\lambda$-harmonique sur $C$, donc $$\lambda_0^N(C)\geq \lambda_0(M_1).$$
\end{proof}

On obtient donc un revêtement riemannien $M\ra M_1$ de groupe de transformation $\Gamma$, sur lequel grace à ce lemme nous allons pouvoir utiliser les méthodes développées à la Section \ref{sec:GInf}. Nous sommes alors très proches des travaux de R. Brooks, pionnier sur les liens entre moyennabilité de groupe de revêtements et spectre du Laplacien. Décrivons ces travaux un peu plus en détail, ainsi que les résultats analogues que nous obtenons.

\paragraph{Bas du spectre d'un revêtement riemannien}
L'ensemble de ce paragraphe est tiré de \cite{Br2}. Soit $M\ra M_1$ un revêtement riemannien de variétés (de dimension finie quelconque), galoisien de groupe de transformation $\Gamma=\pi_1(M_1)/\pi_1(M)$. On suppose que $M_1$ admet une première fonction propre $\psi$ strictement positive, associée au bas du spectre $\lambda_0(M_1)$. On note encore $\psi$ son relevé à $M$, et pour tout domaine fondamental $F$ dans $M$ pour l'action de $\Gamma$ on considère la propriété suivante :\\
\\
\emph{(Br) : Il existe un compact $K\subset F$ tel que 
$$\g h_\psi(F\bs K) = \inf_S\left\{\frac{\int_S \psi^2 d{\rm area}}{\int_{{\rm int}(S)}\psi^2d{\rm vol}}\right\}>0,$$
où $S$ parcours l'ensemble des hypersurfaces découpant $F$ en une partie compacte et une partie non compacte, avec $S\cap K = \vd$, et ${\rm int}(S)$ désigne la composante compacte de $F\bs S$.}\\
\\
Le résultat principal de \cite{Br2} est alors

\begin{theo}[Brooks, 86, Thm 2]
Sous les hypothèses précédentes, si $M_1$ possède un domaine fondamental $F$ pour l'action de $\Gamma$ vérifiant la propriété \emph{(Br)}, alors $$\lambda_0(M_1) \geq \lambda_0(M_2)$$ avec égalité si et seulement $\Gamma$ est moyennable.
\end{theo}

L'hypothèse \emph{(Br)}, peu explicite, implique entre autres que le trou spectral de $M_1$ est positif. Brooks conjecture que ce serait là une hypothèse suffisante pour obtenir ce résultat. Il étudie alors une situation où il est en mesure de contrôler cette hypothèse : il montre que si $M_1$ est une variété hyperbolique géométriquement finie sans cusp à trou spectral positif, alors il existe un domaine fondamental dans $M$ vérifiant l'hypothèse \emph{(Br)}. Il obtient donc dans ce cas :

\begin{theo}[Brooks, 86, Thm 3]
Soit $M_1$ une variété hyperbolique géométriquement finie sans cusp à trou spectral positif, si $M\ra M_1$ est un revêtement riemannien galoisien, alors $\lambda_0(M_1) \geq \lambda_0(M_2)$ avec égalité si et seulement si $\pi_1(M_1)/\pi_1(M)$ est moyennable.
\end{theo}

Dans le cas des surfaces hyperboliques, notons quelques différences entre l'approche de Brooks et la nôtre. Tout d'abord, nous travaillons sur le bas du spectre d'un domaine fondamental à bord géodésique avec conditions de Neumann, tandis que Brooks regarde celui du quotient. La généralisation suivante du Lemme \ref{lemm:SpecDomFond1} lève en partie cette différence :

\begin{lemm}
Soit $M\ra M_1$ un revêtement riemannien galoisien de groupe $\Gamma$ et $F$ un domaine fondamental dans $M$ pour l'action de $\Gamma$, on suppose que $F$ est connexe et $\Cl C^1$ par morceaux. Alors 
$$\lambda_0^N(F) = \lambda_0(M_1).$$
\end{lemm}
\begin{proof}
Il suffit de recopier, mot pour mot, la démonstration du Lemme \ref{lemm:SpecDomFond1}.
\end{proof}

Ensuite, nous travaillons sur des surfaces modelées sur des graphes qui ne sont pas nécessairement les graphes de Cayley de groupes de type fini : cela semble plus général que le résultat de Brooks. En fait, la méthode qu'il utilise s'adapte très bien aux cas que nous traitons, \emph{à condition que la cellule n'ait pas de cusp}.

Une limite de notre méthode semble venir de ce que nous utilisons dans notre construction l'invariance de la cellule par une isométrie qui échange les composantes de bord, ce dont n'a pas besoin Brooks. Mais comme nous l'avons noté à la Section \ref{ssec:GInfBorn}, cette hypothèse ne nous sert qu'à recoller les fonctions propres sur les bords des cellules : dans le cas d'un revêtement, elle est \emph{superflue}. En effet, si $\psi$ est la première fonction propre de $M_1$, lorsqu'on la relève en $\tilde \psi$ sur $M$, on obtient une fonction $\psi$ sur le domaine fondamental $F$ (qui sera notre cellule $C$) qui se recolle évidemment d'une cellule à l'autre puisqu'elle est continue sur $M$. Notre méthode s'adapte alors sans difficulté pour montrer le Théorème \ref{th:GInfRev} :

\begin{theo}
Soit $M_1$ une surface hyperbolique à trou spectral positif, et $M\ra M_1$ un revêtement riemannien galoisien de groupe de revêtement $\Gamma$ de type fini. Supposons qu'il existe un domaine fondamental $F$ dans $M$ pour l'action de $\Gamma$ dont le bord est une union de géodésiques fermées.

Alors il existe des constantes $A_1$ et $A_2$ ne dépendant que de propriétés spectrales de $M_1$ et de la longueur des composantes de $\bd F$ telles que
$$\lambda_0(M_1)+A_1\mu_0(\Gamma)\leq \lambda_0(M)\leq \lambda_0(M_1)+A_2 \g h(\Gamma).$$
\end{theo}

On a noté $\g h(\Gamma)$ et $\mu_0(\Gamma)$ la constante de Cheeger et le bas du spectre du graphe de Cayley associé à un système fini quelconque de générateurs de $\Gamma$. 

\begin{rema}
Ce résultat présente deux améliorations nettes par rapport au résultat de Brooks : notre méthode donne un \emph{contrôle explicite de} $\lambda_0(M)$ en fonction des constantes du groupe $\Gamma$ et de la variété $M_1$, et nous autorisons la \emph{présence de cusps}. En effet, on montre que l'hypothèse \emph{(Br)} n'est pas nécessairement vérifiée en présence de cusp. 
\end{rema}

\begin{rema}\label{rem:DomPolyh}
L'hypothèse que $\bd F$ \emph{est totalement géodésique} est, elle, une vraie limitation de la portée de notre résultat. En particulier, c'est elle qui empèche la généralisation immédiate de nos méthodes en dimension supérieure : lorsque $\Gamma$ est un groupe d'isométrie de $\Bb H^n$ de type fini, on peut souvent se ramener à un domaine fondamental dont le bord est \emph{polyédral} à faces et arêtes totalement géodésiques, mais on ne peut supposer en général que ce bord est totalement géodésique. Il est naturel de penser que nos méthodes s'adaptent à cette situation, y compris lorsque la courbure n'est que négative pincée (et non constante), mais cela implique de contrôler de près ce qu'il se passe aux \emph{angles} du bord du domaine fondamental. Cela fera l'objet d'une prochaine étude.
\end{rema}

\subsection{Perspectives}

La Remarque \ref{rem:DomPolyh} ci-dessus présente la première généralisation de notre méthode qu'il semble naturel de mener. Nous présentons maintenant deux autres directions pour poursuivre cette étude ; la première fait suite à notre remarque et s'intéresse à des variétés de dimension quelconque, tandis que la seconde cherche à affiner nos résultats sur certaines surfaces de genre infini.

\paragraph{Découpages bornés généraux}

On peut se demander à quel point le Théorème \ref{th:GInfBorn} peut se généraliser pour caractériser la nullité du bas du spectre des variétés de volume infini. Certaines variétés riemanniennes, pas forcément de dimension 2 ni à courbure constante, admettent un découpage en \emph{polyèdres} dont le nombre de côtés est borné (par exemple une triangulation), dont les faces sont totalement géodésiques, dont les volumes des $n$-simplexes et des $(n-1)$-faces sont bornés, et dont les trous spectraux sont uniformément bornés. Nous continuerons à appeler cela un \emph{découpage borné} de la variété, et le graphe sous-jacent se définit de la même façon que précédemment. La question naturelle est la suivante :

\begin{quest}
Soit $M$ une variété riemannienne admettant un découpage borné, a-t-on $\lambda_0(M) = 0$ si et seulement si le graphe sous-jacent au découpage est moyennable ?
\end{quest}

Une adaptation des méthodes précédentes donnera peut-être une réponse affirmative ; elle nécessiterait ainsi qu'il l'a été noté à la Remarque \ref{rem:DomPolyh} de s'intéresser à la validité de nos méthodes lorsque nos cellules sont des domaines à bord géodésique par morceau. Dans le cas d'une réponse affirmative, il serait alors intéressant d'obtenir une classe de variétés, la plus large possible, qui admettent un découpage borné : nous obtiendrions pour ces variétés une caractérisation combinatoire naturelle de la nullité du spectre du Laplacien.

\begin{quest}
Quelles variétés riemanniennes admettent un découpage borné ?
\end{quest}

\paragraph{Pincement d'un nombre infini de géodésiques}

Enfin, à l'aide de nos méthodes, on peut espérer adapter les résultats de [Colbois] et [Colbois-Colin] rappelés à la Section \ref{ssec:PincGeod} à certaines surfaces de genre infini. En particulier, dans le cas d'une surface $M$ admettant un découpage contrôlé (voir Section \ref{ssec:GInfBorn}) que l'on pince uniformément le long des géodésiques de découpages, on cherchera à réécrire le Théorème \ref{th:EqCol}  pour le bas du spectre des surfaces $M_\epsilon$ obtenues par ce pincement. Dans le cas d'une surface modelée sur un graphe ou d'un revêtement lorsque le volume des cellules est infini, il sera intéressant de chercher ce que deviennent les théorèmes \ref{th:GInf2} et \ref{th:GInfRev} lorsque l'on pince la surface : on peut alors espérer obtenir un équivalent de 
$$\delta_\epsilon = \lambda_0(M_\epsilon)-\lambda_0^N(C_\epsilon)$$
lorsque $\epsilon\ra 0$. On a noté ici $M_\epsilon$ et $C_\epsilon$ les surfaces obtenues en pinçant uniformément $M$ et $C$ le long des bords des cellules.

\appendix
\section{Appendice : bas du spectre avec condition de Neumann au bord}\label{append}

\subsection{Caractérisation par le spectre positif}

Soit $M$ une variété complète non compacte, dont le bord $\bd M$ est compact et $\Cl C^1$ par morceaux. On rappelle que nous disons qu'une fonction $f$ sur $M$ à valeur réelle vérifie les \emph{conditions de Neumann sur} $M$ si et seulement si pour tout $\xi\in\bd M$, 
$$\frac{\bd f}{\bd \nu}(\xi) = g_\xi(\nabla f(\xi),\nu(\xi)) = 0$$ où $\nu(\xi)$ est la normale à $\bd M$ en $\xi$. Notre objectif est de démontrer le Théorème \ref{th:CarSul}, que nous adaptons de \cite{Sul} et qui donne une caractérisation importante du bas du spectre de $M$ avec condition de Neumann :

\begin{theo}
Pour tout réel $\lambda$, il existe une fonction $\phi$ $\Cl C^\infty$ \emph{$\lambda$-harmonique positive} sur $M$ avec condition de Neumann sur $\bd M$ si et seulement si $\lambda\leq\lambda_0^N(M)$.
\end{theo}
On rappelle que nous notons
$$\lambda_0^N(M) = \inf_f\frac{||\nabla f||^2_{L^2(M)}}{||f||_{^2L^2(M)}}$$
où $f$ parcourt l'ensemble des fonctions $\Cl C^\infty$ à support compact dans $M$. On peut donc réécrire ce théorème sous la forme

$$\lambda_0^N(M) = \sup\{\lambda\in\Bb R : \exists f\in\Cl C^\infty(C), f>0 : \Delta f = \lambda f\}.$$

Lorsque $\bd M = \vd$, ce résultat est exactement le Théorème 2.1 de \cite{Sul}. La démonstration que nous présentons est adaptée de celle de \cite{Sul}, §3-4, bien connue quoique fort peu détaillée dans cet article. Elle utilise ce que l'on appelle communément le \emph{mouvement brownien} : ce terme vient de ce que, dans un modèle physique statistique du type de celui du gaz parfait, $p(x,y,t)dV(y)$ est la densité de probabilité, pour une particule qui se trouvait en $x$ à $t=0$, de se trouver au voisinage de $y$ au temps $t$.

Soit $K$ un voisinage compact de $\bd M$ dans $M$, et $(M_j)_{j\in\Bb N}$ une famille croissante d'ouverts relatifs de $M$ contenant $K$, d'adhérence compacte et tels que
$$\bigcup_j M_j = M.$$ On note $\bd^1 M_j = \bd M\subset M_j$ et $\bd^2 M_j = \bd M_j\bs\bd M$ qu'on suppose également $\Cl C^1$ par morceaux.
On notera encore $$\lambda_0^j=\lambda_0(M_j) = \inf_f\frac{||\nabla f||^2_{L^2(M_j)}}{||f||^2_{L^2(M_j)}}$$
où $f$ parcourt l'ensemble des fonctions $\Cl C^\infty$ à support compact dans $M_j=\bar{M}_j\bs \bd^2 M_j$, cela correspond au bas du spectre du Laplacien avec \emph{condition de Neumann} sur $\bd^1 M_j$ et \emph{de Dirichlet} sur $\bd^2 M_j$.
On a alors $$\lambda_0(M) = \inf_j\lambda_0^j.$$

\pgh

On note $p^j(x,y,t) = p_{M_j}(x,y,t)$ le noyau de la chaleur de $M_j$ associé au problème mixte considéré, c'est à dire la solution fondamentale de l'équation aux dérivées partielles
$$\Delta f = -\frac{\bd f}{\bd t}.$$
On rappelle que nous utilisons un Laplacien défini positif, qui s'écrit sur $\Bb R^n$ en coordonnées euclidiennes
$$\Delta = -\sum_{i=1}^n\frac{\bd^2}{\bd x_i^2}.$$
On a alors pour tous $x,y\in M_j$ et $t>0$,
\beq \label{eq:NoyChal} p^j(x,y,t) = \sum_k e^{-\lambda_k^j t}\phi_k^j(x)\phi_k^j(y),\eeq
où $\phi_k^j$ est la fonction propre du Laplacien avec condition de Neumann sur $\bd^1M_j$ et de Dirichlet sur $\bd^2M_j$ associée à la valeur propre $\lambda_k^j$. Pour tous $x,y\in M$ tels que $x,y\in M_j$ pour tous $j\geq j_0$, on note
$$p(x,y,t) = \inf_{j\geq j_0} p^j(x,y,t).$$
On appelle $p(x,y,t)$ ainsi défini le \emph{noyau de la chaleur minimal} associé au problème de Neumann sur $M$.

\subsection{Construction de fonctions $\lambda$-harmoniques par diffusion}

Les énoncés donnés dans ce paragraphe nécessiteraient, pour être prouvés, un développement de la \emph{théorie de la diffusion} associée à un opérateur elliptique qui passe entre autres par les \emph{intégrales stochastiques} bien plus long que ce qui est souhaitable ici.  Nous nous contentons donc de présenter certaines définitions, les résultats et les idées clés de la démonstration. Nous invitons le lecteur à se référer à \cite{Sul} pour une présentation analogue à la nôtre dans le cas du problème sans bord, à \cite{Mal} et \cite{Cha} pour la construction du mouvement brownien à l'aide du noyau de la chaleur. Les bases de probabilités nécessaires à cette démonstration se trouvent par exemple dans \cite{Ba1}, chapître I, et le détail de nos démonstrations à partir des intégrales stochastiques se trouve dans \cite{Ba2} dans le cas d'ouverts de $\Bb R^d$. La justification de leur adaptation aux variétés riemanniennes se trouve par exemple dans \cite{Em}. 

\begin{defi}
Soit $j>0$ et $M_j$ l'un domaine de $M$ défini ci-dessus. Notons $\Omega$ l'ensemble des chemins continus de $\bar{\Bb R}_+$ dans $M$ et $\tau : \Omega\ra \bar{\Bb R}_+$ défini pour tout $\omega\in\Omega$ par
\beq \label{eq:Tau} \tau(\omega) = inf\{t>0 : \omega(t)\in\bd^2 M_j\}.\eeq
Soit $\Omega^j$ l'ensemble des chemins de $\bar{\Bb R}_+$ dans $\bar M_j$ tels que 
$$\forall t \geq\tau(\omega), \omega(t)=\omega(\tau(\omega))\in \bd^2 M_j$$
et
$$\Omega^j_x = \{\omega\in \Omega^j :\omega(0) = x\}.$$
On appelle $\Omega^j$ l'ensemble des \emph{trajectoires} dans $M_j$.
\end{defi}

On appelle \emph{cylindre} de $\Omega^j_x$ un ensemble de la forme
$$A = \{\omega\in \Omega^j_x : (\omega(t_1),\ldots,\omega(t_k))\in B\},$$
où $k\in \Bb N$, et $B\subset (M_j)^k$ est un borélien, et les $t_j$ sont des réels
$$0\leq t_1<t_2...<t_k.$$

Pour tout cylindre $A$ de la forme précédente, on pose
$$\Bb P^j_x(A) = \int_Bp^j(x,y_1,t_1)p^j(y_1,y_2,t_2-t_1)...p^j(y_{k-1},y_k,t_k-t_{k-1})dV(y_1)...dV(y_k),$$
où $dV$ désigne la mesure canonique associée à la métrique de $M$. On peut montrer à partir de la propriété de semi-groupe des noyaux de la chaleur que $\Bb P^j_x$ s'étend en une unique \emph{mesure de probabilité} sur la $\sigma$-algèbre de $\Omega^j_x$ engendrée par ses cylindres (voir \cite{Cha}). $\Omega^j$ est \emph{l'espace de Wiener} sur $M_j$, et $\Bb P_x^j$ \emph{la mesure de Wiener} en $x$. 

On considère le processus aléatoire $(X_t)_{t\geq0}$ sur $(\Omega_x^j,\Bb P^j_x)$ défini pour tout $\omega\in\Omega_x^j$ par
$$X_t(\omega) = \omega(t).$$
D'après la définition de $\Bb P_x^j$, on a pour tout borélien $B$ de $M$
$$\Bb P_x^j(X_{t+s}\in B|X_t = z) = \int_Bp^j(z,y,s)dV(y) = \Bb P_z^j(\{X_s\in B\})\ :$$
$X_t$ est un processus de Markov de loi $p^j$.

On appelle \emph{mouvement brownien sur $M_j$} (avec réflexion sur $\bd^1 M_j$, ce qui sera désormais sous-entendu) le processus aléatoire $(X_t)_{t\geq 0}$ muni de la loi $\Bb P_x^j$. 

Soit $f : M^j\ra \Bb R$ une fonction de classe $\Cl C^2$. En notant pour tous vecteurs $Y,Z\in T_xM$, $g_x(Y,Z) = Y.Z$, la formule d'Itô pour le mouvement Brownien $(X_t)_{t\geq 0}$ s'écrit (voir \cite{Ba1} p 49, \cite{Em} p 34) s'écrit :

\beq \label{eq:Ito1} f(X_t) = f(X_0) + \int_0^t \nabla f(X_s).dX_s-\int_0^t\Delta f(X_s)ds.\eeq

Le dernier terme de notre formule diffère de la formule de \cite{Ba1} p 49 d'un facteur $-2$ : cela vient de ce que notre convention de signe pour le Laplacien est opposée à celle de Bass, et de ce que le mouvement brownien habituellement considéré lors de l'écriture de la formule d'Itô a pour probabilités de transitions la solution élémentaire de l'équation
$$-\frac{1}{2}\Delta f = \frac{\bd f}{\bd t},$$
alors que nous ne gardons pas ce facteur $\frac{1}{2}$ dans notre construction (voir par exemple \cite{Ba2} p53).

D'après \cite{Ba2} p 33, comme $X_t$ est un brownien avec reflexion normale sur $\bd^1 M_j$, on peut écrire
$$dX_t = dW_t + \nu(X_t)dL_t,$$
où $W_t$ est un brownien sans reflexion sur $M_j$, $\nu(X_t)$ est la normale rentrante à $\bd^1M_j$ en $X_t$ lorsque $X_t\in\bd^1M_j$, et $0$ ailleurs, et $L_t$ le \emph{temps local} sur $\bd^1M_j$. Ce temps local est un processus positif croissant à variation bornée, strictement croissant lorsque $X_t\in\bd^1M$, défini par
$$L_t = \lim_{\epsilon\ra 0} \frac{1}{\epsilon}\int_0^t \textbf{1}_{d(X_s,\bd^1M_j)\leq\epsilon}ds,$$
où $d(X_s,\bd^1M_j)$ désigne la distance de $X_s$ à $\bd^1M_j$ pour la distance induite par la métrique sur $M$. La formule (\ref{eq:Ito1}) devient alors

\beq \label{eq:Ito2} f(X_t) = f(X_0)+\int_0^t \nabla f(X_s).dW_s+\int_0^t\nabla f(X_t).\nu(X_t)dL_t-\int_0^t\Delta f(X_s)ds,\eeq

Supposons $f$ harmonique sur $M^j$, avec condition de Neumann sur $\bd^1 M^j$, c'est à dire
$$\nabla f.\nu \equiv 0$$
sur $\bd^1M^j$. En intégrant (\ref{eq:Ito2}) sur $\Omega_x^j$, on obtient

\begin{theo}
Pour toute fonction $f$ harmonique sur $M^j$ avec condition de Neumann sur $\bd^1 M^j$ et pour tout $t>0$, 
$$ f(x) = \Bb E_x^j(f(X_t)) = \int_{\Omega_x^j}f(\omega(t))d\Bb P^j_x(\omega).$$
\end{theo}

\begin{proof}
Il suffit de remarquer que $$\Bb E_x^j(f(X_0)) = f(x)$$ par définition du mouvement brownien issu de $x$ et que $$\int_0^t \nabla f(X_s).dW_s$$ est une martingale nulle en $t=0$. Les deux autres termes de (\ref{eq:Ito2}) disparaissent pour $f$ harmonique sur $M^j$ avec condition de Neumann sur $\bd^1 M^j$.
\end{proof}

Considérons la variable aléatoire $\tau : \Omega\ra \bar{\Bb R}_+$ définie en (\ref{eq:Tau}). C'est un \emph{temps d'arrêt} (voir \cite{Ba1} p 13) vérifiant $$\Bb P_x^j(\tau>t) = \int_M p^j(x,y,t)dV(y).$$
Puisque 
$$p^j(x,y,t) = \sum_k e^{-\lambda_k^j t}\phi_k^j(x)\phi_k^j(y),$$ on a
$$\lim_{t\ra\infty} e^{\lambda_0^j t}p^j(x,y,t) = \phi_0^j(x)\phi_0^j(y).$$
En particulier, comme $\phi_0>0$ sur $\inter{M}_j$,
\beq \label{eq:TFinPS} \Bb P^j_x(\{\tau> t\}) \sim C e^{-\lambda_0^j t} : \eeq
$\tau$ est fini presque sûrement. Le théorème d'arrêt de Doob (voir \cite{Ba1} p 29) nous donne alors

\begin{theo}\label{th:Harm2}
Pour toute fonction $f$ harmonique sur $M^j$ avec condition de Neumann sur $\bd^1 M^j$,
$$f(x) = \Bb E_x^j(f(X_\tau)) = \int_{\Omega_x^j}f(\omega(\tau(\omega)))d\Bb P_x^j(\omega) = \int_{\bd M_j}f(\xi)d\mu_{j,x}(\xi),$$
\end{theo}

où $\mu_{j,x}$ est la mesure de probabilité sur $\bd^2 M_j$ définie par
$$\mu_{j,x}(B) = \Bb E_x^j(\textbf{1}_{\{X_\tau\in B\}}) = \Bb P_x^j\{\omega\in \Omega_x^j | \exists t>0, \omega(t)\in B\}$$
pour tout borélien $B$ de $\bd M_j$.

Si $f$ est positive et non identiquement nulle sur $\bd^2 M_j$, on peut la prolonger à l'aide du Théorème \ref{th:Harm2} en une fonction harmonique que nous notons encore $f$ sur $M_j$ avec conditions de Neumann sur $\bd^1 M_j$ et continue sur $\bar M_j$. On montre alors montre que $f$ est strictement positive sur l'intérieur de $M_j$. On appelle $\mu_{j,x}$ la \emph{mesure de Poisson}, ou \emph{mesure harmonique}, sur $\bd^2 M_j$ (avec réflexion sur $\bd^1 M_j$) issue de $x$. La construction que nous venons d'en présenter s'appelle la \emph{méthode de balayage de Poincaré}.

\pgh

Soit $\lambda<\lambda_0(M_j)$, nous modifions légèrement cette démonstration pour obtenir des fonctions $\lambda$-harmoniques sur $M_j$. On considère le processus aléatoire $(Y_t)_{t\geq0}$ sur $\Omega_x^j$ défini par
$$Y_t = e^{\lambda t} f(X_t),$$
où $X_t$ est toujours le mouvement brownien avec reflexion défini précédemment. La formule d'Itô s'écrit désormais
\beq \label{eq:Ito3} e^{\lambda t} f(X_t) = f(X_0) +\int_0^t e^{\lambda s} \nabla f(X_s).dX_s+ \int_0^t \lambda e^{\lambda s} f(X_s)ds -\int_0^te^{\lambda s}\Delta f(X_s)ds.\eeq

De même que précédemment, si $f$ est $\lambda$-harmonique (i.e. $\Delta f = \lambda f$) avec condition de Neumann en $\bd^1 M^j$, et en intégrant (\ref{eq:Ito3}) sur $\Omega_x^j$ on obtient

\begin{theo}\label{th:LHarm}
Pour toute fonction $\lambda$-harmonique $f$ sur $M^j$ avec condition de Neumann sur $\bd^1 M^j$ et pour tout $t>0$, 
$$f(x) = \Bb E_x^j(e^{\lambda t}f(X_t)) = \int_{\Omega_x^j}e^{\lambda t}f(\omega(t))d\Bb P_x^j(\omega).$$
\end{theo}

Pour tout $\lambda<\lambda_0^j$, d'après (\ref{eq:TFinPS}) $e^{\lambda \tau}f(X_\tau)$ est sommable : on obtient donc de même que précédemment :

\begin{theo}\label{th:LHarm2}
Pour tout $\lambda<\lambda_0^j$, et pour toute fonction $\lambda$-harmonique $f$ sur $M^j$ avec condition de Neumann en $\bd^1 M^j$,
$$f(x) = \Bb E_x^j(e^{\lambda \tau}f(X_\tau)) = \int_{\Omega_x^j}e^{\lambda\tau(\omega)}f(\omega(\tau(\omega)))d\Bb P_x^j(\omega) = \int_{\bd M_j}f(\xi)d\mu^\lambda_{j,x}(\xi),$$
\end{theo}
où $\mu_{j,x}^\lambda$ est la mesure (finie, non normalisée) définie pour tout borélien $B$ de $\bd M_j$ par :
$$\mu_{j,x}^\lambda(B) = \Bb E_x^j(e^{\lambda\tau} \textbf{1}_{\{X_\tau\in B\}})= \int_{\Omega_x^j} e^{\lambda \tau(\omega)} \textbf{1}_{\{\omega(\tau(\omega))\in B\}} d\Bb P_x^j(\omega).$$

Pour tout $j$, pour tout $\lambda<\lambda_0(M)\leq\lambda_0(M_j)$, soit $f_j$ une fonction positive non identiquement nulle sur $\bd^2 M_j$. Grace au Théorème \ref{th:LHarm2}, on peut la prolonger en une fonction $\lambda$-harmonique strictement positive que nous noterons toujours $f_j$ sur l'intérieur de $M_j$, continue sur $\overline{M_j}$, avec condition de Neumann sur $\bd^1 M_j = \bd M$. Il s'agit maintenant de faire converger une suite de telles fonctions $(f_j)_{j\geq0}$ vers une fonction $\lambda$-harmonique (avec condition de Neumann) définie globalement sur $M$ : c'est l'objet du prochain paragraphe.

\subsection{Principe de Harnack et démonstration du Théorème de Sullivan}

Fixons $x_0\in K\subset M_j$ pour tout $j$. Le résultat suivant, que nous ne démontrerons pas, est tiré de \cite{Sul}, p 336 lorsque $\bd^1M = \vd$. Son analogue avec condition de Neumann sur $\bd^1M$ se montre exactement de la même façon.

\begin{prop}
Soit $x_0\in M_j$ fixé, et $x\in M_j$.

Les mesures $\mu_{j,x}$ et $\mu_{j,x_0}$ sont équivalentes, et la dérivée de Radon-Nikodym
$$\psi_{j,x}(\xi) = \frac{d\mu_{j,x}}{d\mu_{j,x_0}}(\xi)$$
est telle que :

à $\xi$ fixé, la fonction $x\mapsto\psi_{j,x}(\xi)$ est positive et harmonique sur $M_j$, et se prolonge en une fonction continue sur $\overline{M_j}\bs\{\xi\}$, nulle sur $\bd^2 M_j\bs\{\xi\}$ et qui admet un pôle en $\xi$.

De même, si on pose $$\psi_{j,x}^\lambda(\xi) = \frac{d\mu_{j,x}^\lambda}{d\mu_{j,x_0}^\lambda}(\xi),$$
alors pour tout $\xi$ fixé la fonction $x\mapsto\psi_{j,x}(\xi)^\lambda$ est positive et $\lambda$-harmonique sur $M_j$, et se prolonge en une fonction continue sur $\overline{M_j}\bs\{\xi\}$, nulle sur $\bd^2 M_j\bs\{\xi\}$ et qui admet un pôle en $\xi$.
\end{prop}

\begin{coro}[Principe de Harnack]
Soit $\lambda<\lambda_0^j$, soit $f$ une fonction positive $\lambda$-harmonique sur $M_j$, continue sur $\overline{M_j}$. On a alors
$$f(x) = \int_{\bd M_j}f(\xi)\psi_{j,x}^\lambda(\xi)d\mu_{j,x_0}^\lambda(\xi).$$
En particulier, dans un compact de $M_j$ contenant $x_0$, les valeurs de $f$ sont des combinaisons convexes de $\psi_{j,x}^\lambda(\xi)$ à coefficients en $\xi$ fixes, et $\psi_{j,x}^\lambda(\xi)$ est bornée (en $x$) sur ce compact.
\end{coro}

Ce principe de Harnack va nous permettre d'obtenir une convergence uniforme sur tout compact d'une suite de fonctions $\lambda$-harmoniques.

\begin{proof}[Démonstration du Théorème de Sullivan]

On sait que la suite $\lambda_0^j = \lambda_0(M_j)$ est décroissante, et qu'on a
$$\lambda_0(M) = \inf_j \lambda_0(M_j).$$

Soit $i>0$, et $\lambda<\lambda_0$, on a $\lambda<\lambda_0^j$ pour tout $j>0$. D'après la construction précédente, pour tout $j>i$, il existe une fonction $\lambda$-harmonique $f_j$ strictement positive sur $M_j$ avec condition de Neumann sur $\bd M$ et valant $1$ en $x_0$. On prolonge les $(f_j)_{j>0}$ en des fonctions continues bornées sur $M$. Comme $M$ est complète, il existe une fonction $f$ vers laquelle (quitte à extraire une sous-suite) la suite des $(f_j)$ converge simplement. Comme les $(f_j)_{j>i}$ sont $\lambda$-harmoniques sur $M_i$, d'après le principe de Harnack ci-dessus la suite converge uniformément vers $f$ sur tout compact de $M_j$ contenant $x_0$. Par convergence dominée, $f$ est $\lambda$-harmonique sur tout compact de $M_j$, avec condition de Neumann sur $\bd M$, ce pour tout $j>0$. 

Pour $\lambda = \lambda_0$, on considère une suite croissante de réels $(\lambda^n)_{n>0}$, avec 
$$\lim_{n\ra\infty}\lambda^n = \lambda_0,$$ et une suite $(f_n)_{n>0}$ de fonctions $\lambda^n$-harmoniques positives sur $M$, valant $1$ en $x_0$, construites comme précédemment. Quitte à extraire une sous-suite, les $(f_n)$ convergent uniformément sur tout compact vers une fonction positive $f_0$ qui, par convergence dominée, est $\lambda_0$-harmonique sur $M$.

\pgh
Réciproquement, soit $\lambda\in \Bb R$ tel qu'il existe  une fonction $\lambda$-harmonique positive $f$ sur $M_j$ (avec condition de Neumann sur $\bd M = \bd^1 M_j$), on a

\begin{lemm}
Pour tout $t>0$,
$$f(x) = \int_{M_j} e^{\lambda t}p^j(x,y,t)f(y)dV(y) + \int_{\Omega^j_x} e^{\lambda\tau(\omega)}f(\omega(\tau(\omega)))\textbf{1}_{\{\tau\leq t\}}d\Bb P_x^j(\omega).$$
\end{lemm}

\begin{proof}
C'est simplement le Théorème \ref{th:LHarm} écrit à l'instant $t>0$.
\end{proof}

On a donc pour tout $t$ positif,
$$f(x)\geq \int_{M_j} e^{\lambda t}p^j(x,y,t)f(y)dV(y).$$
Comme $f(x)$ est finie et $$\lim_{t\ra\infty} e^{\lambda_0^j t}p^j(x,y,t) = \phi_0^j(x)\phi_0^j(y),$$ avec $\phi_0^j>0$ sur $\inter{M_j}$, on a nécessairement $\lambda\leq\lambda_0^j$. Ceci est valable pour tout $j\in\Bb N$, ce qui conclut la preuve du théorème.

\end{proof}

\end{document}